\documentclass[reqno,10pt]{amsart}
\usepackage{amsmath,amssymb}

\newcommand{\twocase}[5]{#1 \begin{cases} #2 & \text{\rm #3}\\ #4
&{\rm #5} \end{cases}   }

\newcommand{\weight}{\frac{\zeta(2)}{L(1,\sym^2f)}}

\newcommand{\un}{\rm U}
\newcommand{\sy}{\rm USp}
\newcommand{\usp}{\rm USp}
\newcommand{\soe}{\rm SO(even)}
\newcommand{\soo}{\rm SO(odd)}
\newcommand{\so}{\rm O}
\renewcommand{\d}{{\mathrm{d}}} 

\newcommand{\sym}{{\rm sym}}
\newcommand{\hg}[1]{\widehat{g}\left(#1\right)}
\newcommand{\hgi}[1]{\widehat{g_i}\left(#1\right)}
\newcommand{\hgo}[1]{\widehat{g_1}\left(#1\right)}
\newcommand{\hgt}[1]{\widehat{g_2}\left(#1\right)}

\newcommand{\apsf}{a_{\phi\times\sym^2f}}

\newcommand{\afsp}{a_{\phi\times\sym^2f}}

\newcommand{\bfsp}{b_{\phi\times\sym^2f}}
\newcommand{\glp}{\lambda_\phi}
\newcommand{\glf}{\lambda_f}

\newcommand\be{\begin{equation}}
\newcommand\ee{\end{equation}}
\newcommand\bp{\begin{proof}}
\newcommand\ep{\end{proof}}

\newcommand\bea{\begin{eqnarray}}
\newcommand\eea{\end{eqnarray}}
\newcommand\bml{\begin{multline}}
\newcommand\eml{\end{multline}}
\newcommand\bal{\begin{align}}
\newcommand\eal{\end{align}}
\newcommand\bi{\begin{itemize}}
\newcommand\ei{\end{itemize}}
\newcommand\ben{\begin{enumerate}}
\newcommand\een{\end{enumerate}}
\newcommand\bc{\begin{center}}
\newcommand\ec{\end{center}}
\newcommand\ba{\begin{array}}
\newcommand\ea{\end{array}}

\newtheorem{thm}{Theorem}[section]

\newtheorem{cor}[thm]{Corollary}
\newtheorem{lem}[thm]{Lemma}

\newtheorem{defi}[thm]{Definition}

\theoremstyle{definition}
\newtheorem{rek}[thm]{Remark}

\newcommand{\R}{\ensuremath{\Bbb{R}}}

\newcommand{\Z}{\ensuremath{\Bbb{Z}}}

\newcommand{\ga}{\alpha}     
\newcommand{\gb}{\beta}      
\newcommand{\gl}{\lambda}

\newcommand{\hphi}{\widehat{\phi}}  

\newcommand{\fof}{\frac{1}{4}}  
\newcommand{\foh}{\frac{1}{2}}  
\newcommand{\fot}{\frac{1}{3}}  

\newcommand{\FD}{\mathcal{F}} 

\newcommand{\bfA}{\mathbb{A}}
\newcommand{\bfQ}{\mathbb{Q}}
\newcommand{\bfC}{\mathbb{C}}
\newcommand{\bfR}{\mathbb{R}}

\newcommand{\kkot}[1]{ \frac{\sin \pi {#1} }{\pi {#1} } }

\newcommand{\hkn}{H_k}

\newcommand{\gltwo}{{\rm GL}(2)}
\newcommand{\glfour}{{\rm GL}(4)}
\newcommand{\glsix}{{\rm GL}(6)}
\newcommand{\gln}{{\rm GL}(n)}
\newcommand{\GL}{{\rm GL}}    
\newcommand{\SL}{{\rm SL}}    

\DeclareMathOperator{\diag}{diag}

\numberwithin{equation}{section}

\begin{document}

\title[The Low Lying Zeros of Two Families of $L$-functions]{The
  Low Lying Zeros of a $\glfour$ and a $\glsix$ family of
  $L$-functions}

\author{Eduardo Due\~{n}ez}
\address{Department of Mathematics, The University of Texas at San Antonio,
San Antonio, TX   78249} \email{eduenez@math.utsa.edu}

\author{Steven J. Miller}
\address{Department of Mathematics, Brown University, 151 Thayer
 Street,
Providence, RI 02912}
 \email{sjmiller@math.brown.edu}



\thanks{We thank Wenzhi Luo and Peter Sarnak for suggesting this
  problem, as well as for many enlightening conversations, and Jim Cogdell for
  valuable comments in writing the Appendix to this article.  The first-named
  author was partly supported by EPSRC Grant N09176.}
\begin{abstract}
  We investigate the large weight ($k\to\infty$) limiting statistics
  for the low lying zeros of a $\glfour$ and a $\glsix$
  family of $L$-functions, $ \{L(s,\phi \times f): f \in \hkn\}$ and $
  \{L(s,\phi \times \sym^2 f): f \in \hkn\}$; here $\phi$ is a fixed even
  Hecke-Maass cusp form and $\hkn$ is a Hecke eigenbasis for the
  space $H_k$ of holomorphic cusp forms of weight $k$ for the full
  modular group. Katz and Sarnak conjecture that the behavior of
  zeros   near the central point should be well modeled by the behavior
  of eigenvalues near $1$ of a classical compact group. By studying the 1- and 2-level densities, we find
  evidence of underlying symplectic and $\soe$ symmetry,
  respectively. This should be contrasted with previous results of
  Iwaniec-Luo-Sarnak for the families $\{L(s,f): f\in \hkn\}$ and
  $\{L(s,\sym^2f): f\in \hkn\}$, where they find evidence of orthogonal and
  symplectic symmetry, respectively. The present examples suggest a
  relation between the symmetry type of a family and that of its
  twistings, which will be further studied in a subsequent
  paper. Both the $\glfour$ and the $\glsix$ families above have all even
  functional equations, and neither is naturally split from an orthogonal
  family. A folklore conjecture states that such families must be
  symplectic, which is true for the first family but false for the second. Thus the
  theory of low lying zeros is more than just a theory of signs of
  functional equations. An analysis of these families suggest that it is the second
  moment of the Satake parameters that determines the symmetry
  group.
\end{abstract}

\maketitle

\tableofcontents



\section{Introduction}

Assuming GRH, the non-trivial zeros of any $L$-function lie on its
critical line, and therefore it is possible to investigate the
statistics of its normalized zeros. The general philosophy, born
out of many examples and proven cases in function fields
\cite{CFKRS,KS1,KS2,KeSn,ILS}, is that the statistical behavior of
eigenvalues of random matrices (resp., random matrix ensembles) is
similar to that of the critical zeros of $L$-functions (resp.,
families of $L$-functions).


The global $n$-level correlations of high zeros of primitive
automorphic cuspidal $L$-functions, assuming a certain technical
restriction, have been found to agree with the corresponding
statistics of the eigenvalues of complex hermitian matrices (the
Gaussian Unitary Ensemble, or GUE)~\cite{Mon,Hej,RS}. If the
technical restriction mentioned above were to be removed, the
results on $n$-level correlations would imply that the
distributions of the normalized neighbor spacings between
consecutive critical zeros of an $L$-function and between GUE
eigenvalues coincide, as has been numerically observed
\cite{Od,Ru1}.  The same correlations describe the global
statistical behavior of the eigenvalues of other matrix ensembles,
most notably of the classical compact groups (orthogonal, unitary,
symplectic). Being insensitive to the effect of finitely many
zeros, these correlations miss the behavior of the low lying
zeros, the zeros near the central point $s=1/2$.

Katz and Sarnak~\cite{KS1,KS2} showed that there is another
statistic that can distinguish between the classical compact groups.
It is the $n$-level density, and it depends only on eigenvalues
near~$1$. In a number of
cases~\cite{FI,Gu,HM,HR2,ILS,Mil2,Ro,Ru2,Yo}, the behavior of the
low lying zeros of families of $L$-functions is found to be in
agreement with that of the eigenvalues near~$1$ for random matrices
in one of the classical compact groups: unitary, symplectic, and
orthogonal (which is further split into $\soe$ and $\soo$). This
correspondence allows us, at least conjecturally, to assign a
definite symmetry type to each family of $L$-functions.

Let $\phi$ be a fixed even Hecke-Maass cusp form and $\hkn$ a Hecke
eigenbasis for the space of holomorphic cusp forms of (even)
weight~$k$ for the full modular group. Iwaniec-Luo-Sarnak \cite{ILS}
proved that as $k \to \infty$, the family $\{f: f\in\hkn\}$ has
$\soe$ or $\soo$ symmetry (depending on whether $k/2$ is even or
odd), and the family $\{\sym^2 f: f\in\hkn\}$ has symplectic
symmetry. We consider the twisted families $\FD_{\phi\times H_k}$
$=$ $\{\phi \times f: f \in \hkn\}$ and $\FD_{\phi\times \sym^2H_k}$
$=$ $\{\phi \times \sym^2 f: f \in \hkn\}$; the family $\{\phi
\times \sym^2 f\}$ arose in the work of Luo-Sarnak \cite{LS}, where
it is shown to be intimately connected with the relation between the
quantum and classical fluctuations of observables on the modular
surface. In both families, all functional equations are even. We
show that the first family has symplectic symmetry, and the second
$\soe$. Explicitly, our main results are

\begin{thm}\label{thm:sympfam} Let $\phi$ be a fixed even Hecke-Maass
cusp form. As $k\to\infty$, for test functions whose Fourier
transform has small but computable support, the $1$-level density
of the family $\FD_{\phi\times H_k}$ only agrees with symplectic
matrices, suggesting that the underlying symmetry of this family
is symplectic (and uniquely so). \end{thm}

\begin{thm}\label{thm:soefam}
Let $\phi$ be a fixed even Hecke-Maass cusp form. As $k\to\infty$,
for test functions $g$ with ${\rm supp}(\widehat{g}) \subset
(-\frac5{24},\frac5{24})$, the $1$-level density of the family
$\FD_{\phi\times\sym^2 H_k}$ only agrees with $\soe$, $\so$ and
$\soo$ matrices. For small but computable support, the $2$-level
density only agrees with $\soe$ matrices, suggesting that the
underlying symmetry of this family is $\soe$ (and uniquely so).
\end{thm}

For families where the signs of the functional equations are all even
and there is no corresponding family with odd functional equations, a
``folklore'' conjecture (for example, see page 2877 of~\cite{KeSn})
states that the symmetry is symplectic, presumably based on the
observation that $\soe$ and $\soo$ symmetries in the examples known to
date arise from splitting orthogonal families according to the sign of
the functional equations. \emph{A priori} the symmetry type of a
family with all functional equations even is either symplectic or
SO(even). All $L$-functions from elements of $\FD_{\phi\times H_k}$
and $\FD_{\phi\times \sym^2H_k}$ have even functional equations, and
neither family seems to naturally arise from splitting sign within a
full orthogonal family.  By calculating the 1-level density we quickly
see the symmetry of the first is symplectic (as predicted); however,
the second family has orthogonal symmetry (we cannot distinguish
between $\soe, \so$ and $\soo$ due to the small-support restriction on
the allowable test functions).  By calculating the 2-level density for
the second family, we can discard $\so$ and $\soo$. Thus our
calculations are only consistent with the symmetry being $\soe$. As
our purpose is to show that the theory of low lying zeros is more than
just a theory of signs of functional equations, we do not concern
ourselves with obtaining optimal bounds in terms of support, instead
simplifying the arguments but still distinguishing the various
classical compact group candidates.

In studying the symmetry groups of $\FD_{\phi\times H_k}$ and
$\FD_{\phi\times\sym^2 H_k}$, we see that twisting a family with
orthogonal (respectively, symplectic) symmetry by a fixed $\gltwo$
form flips the symmetry to symplectic (respectively, orthogonal). The
effect on the symmetry group by $\gln$ twisting (by a fixed form, or
by a second family) in some cases will be described in a subsequent
paper (\cite{DM}). The main result is that, for any family
$\mathcal{F}$ satisfying certain technical conditions, we can attach a
symmetry constant $c_{\mathcal{F}}$, with $c_\mathcal{F} = 0$ $(1,
-1)$ if the family is unitary (symplectic, orthogonal). For such
families $\mathcal{F}$ and $\mathcal{G}$, the family $\mathcal{F}
\times \mathcal{G}$ (Rankin-Selberg convolution) has symmetry constant
$c_{\mathcal{F}\times\mathcal{G}} = c_{\mathcal{F}} \cdot
c_{\mathcal{G}}$ (compare with \cite{KS1}). In other words, the
symmetry of a product of families is the product of the family
symmetries. This is consistent with earlier results and should be
compared, for instance, with Rubinstein's work~\cite{Ru2} on twisting
the symplectic family of quadratic Dirichlet characters by a fixed
${\rm GL}(n)$ form.

We assume the Generalized Riemann Hypothesis for all $L$-functions
encountered. Mostly GRH is used for interpretation purposes (i.e., if
GRH is true than the non-trivial zeros lie on the critical line, and
we may interpret the $n$-level correlations and densities as spacing
statistics between ordered zeros), though in a few places GRH is
assumed to simplify the derivation of needed bounds (though these
bounds can be derived unconditionally at the cost of a more careful
analysis). In \S\ref{sec:prelims} we review the necessary
preliminaries. We concentrate on the more difficult $\glsix$ family in
\S\ref{sec:lphisym2f}, and merely sketch the changes needed to handle
the $\glfour$ family in \S\ref{sec:L-phi-f}; for completeness the
details of the calculation of the gamma factors and signs of the
functional equations are given in Appendix \ref{sec:appendix}. In
\S\ref{sec:summary-all} we analyze our results for these two
families. The evidence suggests that the theory of low lying zeros is
not just a theory of signs of functional equations, but rather more
about the second moment of the Satake parameters. In this regard it is
similar to the universality Rudnick and Sarnak~\cite{RS} found for the
$n$-level correlations of high zeros of a primitive $L$-function
$L(s,\pi)$ ($\pi$ a cuspidal automorphic representation); their
results are a consequence of the universality of the second moments of
the Satake parameters~$a_\pi(p)$.


\section{Preliminaries}\label{sec:prelims}

\subsection{$1$- and $2$-Level Densities}
\label{sec:central-1-2-level} Let $g$ be an even Schwartz test
function on~$\R$ whose Fourier transform
\begin{equation}
  \label{eq:9}
  \widehat{g}(y)\ =\ \int_{-\infty}^\infty g(x) e^{-2\pi ixy}\d x
\end{equation}
has compact support. Let $\mathcal{F}$ be a finite family, all of
whose $L$-functions satisfy GRH. We define the $1$-level density
associated to $\mathcal{F}$ by
\begin{equation}
\label{eq:7} D_{1,\mathcal{F}}(g)\ =\ \frac1{|\mathcal{F}|}
\sum_{f\in\mathcal{F}} \sum_{j} g\left(\frac{\log
c_f}{2\pi}\gamma_f^{(j)}\right),
\end{equation}
where $\foh + i\gamma_f^{(j)}$ runs through the non-trivial zeros of
$L(s,f)$.  Here $c_f$ is the analytic conductor of $f$, and gives
the natural scale for the low zeros. Since $g$ is Schwartz, only low
lying zeros (i.e., zeros within a distance $\ll\frac1{\log c_f}$ of
the central point) contribute significantly. Thus the $1$-level
density is a local statistic which can potentially help identify the
symmetry type of the family.

\begin{rek}\label{rek:weighting1ld}
For technical convenience, as in \cite{ILS,Ro} we will modify
\eqref{eq:7} by weighting each $f$ by a factor $w_f$ which varies
slowly with $f$. These factors simplify applying the Petersson
formula, and (see \cite{ILS}) can be removed at the cost of
additional book-keeping. \end{rek}

Based in part on the function-field analysis where $G(\mathcal{F})$
is the monodromy group associated to the family $\mathcal{F}$, it is
conjectured that for each reasonable irreducible family of
$L$-functions there is an associated symmetry group $G(\mathcal{F})$
(typically one of the following five subgroups of unitary matrices:
unitary~$U$, symplectic~$\usp$, orthogonal~$\so$, $\soe$, $\soo$),
and that the distribution of critical zeros near $\foh$ mirrors the
distribution of eigenvalues near~$1$.  The five groups have
distinguishable $1$-level densities.

To evaluate \eqref{eq:7}, one applies the explicit formula,
converting sums over zeros to sums over primes. Unfortunately, these
prime sums can often only be evaluated for small support. If one
allows test functions with ${\rm supp}(\widehat{g}) \subset
(-\delta,\delta)$, then for any $\delta>0$ the orthogonal,
symplectic and unitary symmetries can be mutually distinguished via
their $1$-level density. However, if $\delta \leq 1$ then the
$1$-level densities of the three orthogonal types O, SO(even)
and~SO(odd) cannot be distinguished from one another.

In order to distinguish between the three orthogonal symmetry
types we study the $2$-level density of the family, defined as
follows.  Let $g(x) = g_1(x_1)g_2(x_2)$, each $\widehat{g_i}$ of
compact support. Then
\begin{equation}
\label{eq:7b} D_{2,\mathcal{F}}(g)\ =\ \frac1{|\FD|}
\sum_{f\in\mathcal{F}} \sum_{j_1 \neq \pm j_2} g_1\left(\frac{\log
c_f}{2\pi}\gamma_f^{(j_1)}\right)g_2\left(\frac{\log
c_f}{2\pi}\gamma_f^{(j_2)}\right),
\end{equation} Miller \cite{Mil1} observed that an
advantage of studying the $2$-level density is that, even for
arbitrarily small support, the three orthogonal types of symmetry
are mutually distinguishable (see \cite{Mil2} where it is used to
discern the symmetry group of families of elliptic curves). An
analogous definition holds for the $n$-level density; as the signs
of our families are constant, our arguments can easily be extended
to determining the $n$-level density (though the support decreases
with $n$).

By \cite{KS1}, the $n$-level densities for the classical compact
groups are
\begin{align}\label{eqdensitykernels}\begin{array}{lcl}
W_{n,\soe}(x) & \ = \ & \det (K_1(x_i,x_j))_{i,j\leq n}\\
W_{n,\soo}(x) & = & \det (K_{-1}(x_i,x_j))_{i,j\leq n} +
\sum_{k=1}^n \delta(x_k) \det(K_{-1}(x_i,x_j))_{i,j\neq k}\\
W_{n,\so}(x) & = & \foh W_{n,\soe}(x) + \foh W_{n,\soo}(x)\\
W_{n,\un}(x) & = & \det (K_0(x_i,x_j))_{i,j\leq n}\\
W_{n,\sy}(x) &=& \det (K_{-1}(x_i,x_j))_{i,j\leq n}, \end{array}
\end{align}
where $K(y) = \kkot{y}$, $ K_\epsilon(x,y) = K(x-y) + \epsilon
K(x+y)$ for $\epsilon = 0, \pm 1$ and $\delta(x)$ is the Dirac delta
functional; see \cite{HM} for a more tractable formula for the
$n$\textsuperscript{th} centered moments for test functions whose
Fourier transforms have support suitably restricted. It is often
more convenient to work with the Fourier transforms of the
densities. For the 1-level densities we have
\begin{align}
\begin{array}{lcl}
\widehat{W}_{1,\soe }(u) &\ =\ & \delta(u) + \foh I(u) \\
\widehat{W}_{1,\soo }(u) & = & \delta(u) - \foh I(u) + 1 \\
\widehat{W}_{1,\so }(u) & = & \delta(u) + \foh \\
\widehat{W}_{1,U }(u) & = & \delta(u) \\
 \widehat{W}_{1,\sy}(u) & = &
\delta(u) - \foh I(u),
\end{array}\end{align}
where $I(u)$ is the characteristic function of $[-1,1]$. The three
orthogonal densities are indistinguishable for test functions of
small support. Explicitly, for test functions $g$ such that ${\rm
supp}(\widehat{g}) \subset (-1,1)$, we have
\begin{align}\label{eq:FTphi1ld}
\begin{array}{lcl}
\int \widehat{g}(u) \widehat{W}_{1,\soe}(u)\d u & = & \widehat{g}(u) + \foh g(0) \\
\int \widehat{g}(u)\widehat{W}_{1,\soo}(u)\d u &\ =\ &
\widehat{g}(u) + \foh g(0) \\ \int
\widehat{g}(u)\widehat{W}_{1,\so}(u)\d u & = & \widehat{g}(u) + \foh
g(0) \\ \int \widehat{g}(u)\widehat{W}_{1,\un}(u)\d u & = &
\widehat{g}(u)
\\ \int \widehat{g}(u)\widehat{W}_{1,\sy}(u)\d u & = & \widehat{g}(u) - \foh
g(0).
 \end{array}\end{align}

We record the effect of the Fourier transform of the 2-level density
kernel on our test functions. Let $c(\mathcal{G}) = 0$ (respectively
$\foh$, $1$) for $\mathcal{G} = \soe$ (respectively $\so$, $\soo$).
For even functions $\widehat{g}_1(u_1)\widehat{g}_2(u_2)$ supported
in $|u_1| + |u_2| < 1$,
\begin{eqnarray}\label{eq:soetwo} \begin{array}{lcl}
\int \int \widehat{g_1}(u_1)\widehat{g_2}(u_2)
\widehat{W_{2,\mathcal{G}}}(u) \d u_1\d u_2 &= &
\Big[\widehat{g}_1(0) + \foh g_1(0) \Big] \Big[\widehat{g}_2(0) +
\foh g_2(0) \Big] \nonumber\\ & & + \ 2 \int |u| \widehat{g}_1(u)
\widehat{g}_2(u)\d u\nonumber\\ & &  -\ 2 \widehat{g_1g_2}(0) -
g_1(0)g_2(0) \nonumber\\ & & + \ c(\mathcal{G})g_1(0)g_2(0).
\end{array}
\end{eqnarray}Thus, for arbitrarily small support, the 2-level
density distinguishes the three orthogonal groups (see \cite{Mil1}
for the calculation).

\subsection{Cusp Forms} \label{sec:spaces-cusp-forms}
We quickly review some facts about cusp forms; see \cite{Iw2,ILS}
for details. Let $S_k$ be the space of holomorphic cusp forms of
weight~$k$ (an even positive integer) and level~$1$ (that is, for
the full modular group $\Gamma=SL(2,\Z)$).  Let $\hkn$ be a basis of
Hecke eigenforms. Then \be \dim S_k\ =\ |\hkn|\ =\ \frac{k}{12} +
O(1). \ee
 Any $f \in H_k$ has a Fourier
expansion \be f(z)\ =\ \sum_{n=1}^\infty a_f(n) e(nz), \ee and we
shall henceforth assume $f$ is normalized so that $a_f(1) = 1$.
Two other useful normalizations for the coefficients are
\begin{align}
\label{eq:35}
\gl_f(n) &\ =\ a_f(n) n^{ -\frac{k-1}{2} } \\
\label{eq:36} \psi_f(n) &\ =\ \sqrt{ \frac{ \Gamma(k-1)}{(4\pi
n)^{k-1} } }\ \frac{1}{\|f\|}\ a_f(n),
\end{align}
with $\|f\|$ the Petersson $L^2$-norm of $f$. As mentioned in Remark
\ref{rek:weighting1ld}, $\psi_f(n)$ will lead to a weighted sum
which simplifies the application of the Petersson formula. Essential
in our investigations will be the multiplicativity properties of the
Fourier coefficients.

\begin{lem} Let $f$ be a cuspidal Hecke eigenform of level~$1$. Then
  \be \gl_f(m) \gl_f(n) \ =\ \sum_{d|(m,n)} \gl_f\left( \frac{mn}{d^2}
  \right). \ee
\end{lem}
In particular, we have
\begin{cor} Let $(m,n)=1$, and $p$ be a prime. Then
  \bea\label{eqlambdaexpand} \gl_f(m) \gl_f(n) &\ =\ & \gl_f(mn) \nonumber\\
  \gl_f(p)^2 & = & \gl_f(p^2) + 1. \eea
\end{cor}

\subsection{Summation Formulas}
\label{sec:petersson-formulae}
We recall some standard formulas for summing Fourier coefficients
over our families and test functions over primes.
\begin{defi}[Diagonal Symbol]
  \bea \Delta_{k}(m,n)  &\ =\ & \sum_{f \in H_k} \psi_f(m)
  \overline{\psi}_f(n)  \nonumber\\ \twocase{\delta(m,n) & = &
  }{1}{if $m = n$}{0}{otherwise.}
  \eea
\end{defi}

We rephrase the results from \cite{ILS} in our language. By their
equations 2.8, 2.52-2.54, and recalling that $|\hkn| = \frac{k}{12}
+ O(1)$, we find \bea \Delta_k(m,n) &  \ = \  &
\frac{\zeta(2)}{|\hkn| + O(1)} \sum_{f\in \hkn}
\frac{\lambda_f(m)\lambda_f(n)}{L(1,\sym^2 f)}. \eea

\begin{lem}[Petersson Formula]\label{lempeterssonfromils}
For $(m,n) = 1$, $m$ and $n$ of at most $b$ factors,
\be\label{eqilstwotwo} \Delta_{k}(m,n)\  =\ \delta(m,n) +
O_b\left( \frac{ m^\fof n^\fof \log mn }{k^\frac{5}{6}} \right).
\ee For $m, n$ as above and $12\pi\sqrt{mn} \leq k$,
\be\label{eqilstwothree} \Delta_k(m,n)\ =\ \delta(m,n) +
O\left(\frac{\sqrt{mn}}{2^k}\right). \ee
\end{lem}

Note \eqref{eqilstwotwo} and \eqref{eqilstwothree} are Corollaries
2.2 and 2.3 of \cite{ILS}.

The Petersson formula allows us to easily evaluate certain weighted
sums of the Fourier coefficients. As $\|f\|$ is related to
$L(1,\sym^2 f)$, the natural weights are the harmonic weights
$\omega_f = \zeta(2)/L(1,\phi\times\sym^2f)$. These weights are
almost constant (see \eqref{eq:52}), and following \cite{ILS} we may
remove these weights in the applications below. See
\S\ref{sec:prelimsgl6case} for more details. We call terms with $m =
n$ \emph{diagonal} terms; the remaining terms are called
\emph{non-diagonal}. For small support, the non-diagonal terms will
not contribute.

\begin{rek} There exist explicit formulas, involving
Bessel functions and Kloosterman sums, for the error terms in the
expansion of $\Delta_k(m,n)$. For the families studied in
\cite{ILS}, by analyzing these terms' contributions they are able
to work with test functions with support greater than $[-1,1]$,
and hence distinguish $\soe$ from $\soo$; see also \cite{HM} where
these terms are handled for the $n$-level densities. Increasing
the support has applications to non-vanishing results at the
central point. We will not be able to obtain such large support
for our families; however, by studying the $2$-level density, we
can still distinguish $\soe$ from $\soo$.
\end{rek}

The following are immediate applications of the Prime Number
Theorem:

\begin{lem}\label{thmprimesums} Let $\widehat{F}$ be an even Schwartz function of
compact support. Then for any positive integer $a$,
\begin{align}
\label{eq:19} \sum_p \widehat{F}\left( a \frac{\log p}{\log R}
\right) \frac{\log p}{\log R} \frac{1}{p} & \ = \ \frac{1}{2a}
F(0) +
O\left( \frac{1}{\log R} \right) \\
\label{eq:23} \sum_p \widehat{F}\left( \frac{\log p}{\log R} \right)
\frac{4 \log^2 p}{\log^2 R} \frac{1}{p} & \ = \ 2
\int_{-\infty}^\infty |u| \widehat{F}(u)\d u + O\left( \frac{1}{\log
R} \right).
\end{align}
\end{lem}

\section{$\mathcal{F}_{\phi\times \sym^2 H_k} = \{\phi\times\sym^2f: f \in H_k\}$}
\label{sec:lphisym2f}

We provide evidence that the underlying symmetry of the family
$\FD_{\phi\times \sym^2 H_k}=\{\phi \times \sym^2 f:f\in H_k\}$ is
$\soe$. In \S\ref{sec:defn-gamma-fe} we calculate the needed
quantities to investigate the distribution of the low lying zeros.
In \S\ref{subseconelevelsix} we calculate the one-level density
for test functions whose Fourier transform has small support,
proving the first half of Theorem \ref{thm:soefam}. Although the
support is contained in $(-1,1)$, the evidence is enough to
discard the possibility of symplectic (or unitary) symmetry;
however the $1$-level density in this range cannot distinguish
between $\so,\soe$ and $\soo$ (even though the even functional
equations suggest, of course, that $\soe$ is the type). To rectify
this deficiency, in \S\ref{sec:2-level-density} we calculate the
2-level density for small support; this suffices to eliminate
$\so$ and $\soo$ and will complete the proof of Theorem
\ref{thm:soefam}.

\subsection{Definition, Gamma Factors, Functional Equation}
\label{sec:defn-gamma-fe}

We use the notation of \S\ref{sec:spaces-cusp-forms} for the
holomorphic Hecke eigenform $f$ and its Hecke eigenvalues
$\lambda_f(n)$.  Let $\phi$ be a fixed even Hecke-Maass cuspidal
eigenform with Laplacian eigenvalue $\lambda_\phi=\frac14+t_\phi^2$
for the full modular group $\Gamma=\text{SL}(2,\Z)$. We normalize
$\phi$ so that $a_\phi(1)=1$, and denote by $\lambda_\phi(n)$ the
corresponding Hecke eigenvalues.

For any unramified prime~$p$, the Satake parameters (of the
principal-series representation of $GL_2(\bfQ_p)$) associated to~$f$
are two complex numbers $\alpha_p,\tilde\alpha_p=\alpha_p^{-1}$
satisfying
\begin{equation}
  \label{eq:0}
  \lambda_f(p^\nu) \ = \ \sum_{\ell=1}^\nu \alpha_p^\ell{\tilde\alpha_p}^{\nu-\ell}.
\end{equation} Since~$f$ is of level~$1$, every prime is unramified.  By the
work of Deligne, $|\alpha_p|=1$ ---the local representation is tempered--- so
that in fact $\tilde\alpha_p$ is the complex conjugate $\overline\alpha_p$ of
$\alpha_p$.  Thus $\lambda_f(p)=\alpha_p+\alpha_p^{-1}$ alone determines
$\alpha_p,\alpha_p^{-1}$. By~\eqref{eq:0}, all of the $\lambda_f(p^\nu)$ are
algebraically expressible in terms of $\lambda_f(p)$ (formula~\eqref{eq:0} is
indeed equivalent to the multiplicativity of the Fourier coefficients).

The Maass form~$\phi$ has Satake parameters
$\beta_p,\tilde\beta_p=\beta_p^{-1}$.  The Ramanujan conjecture states that
$|\beta_p|=1$; while this is still open for Maass forms, powerful bounds
towards Ramanujan are available. Kim and Shahidi \cite{KiSh} proved the
crucial (for us) bound $|\beta_p|,|\tilde\beta_p|\leq p^{\frac5{34}}$.
Observe that $\frac5{34}<\frac16$, which has many important consequences (see
Section~8 of~\cite{KiSh}), and is perhaps not coincidentally all we need
below. The exponent has been recently improved by Kim and Sarnak to
$\frac7{64}$ (see Appendix Two of~\cite{K}).


Denote by $\sym^2f$ be the Gelbart-Jacquet (symmetric-square) lift to (an
automorphic cuspidal representation of) $\GL(3)$ of the cusp
form~$f$ \cite{GeJa}.  Its Fourier coefficients are~\cite{Bu1,Bu2}
\begin{equation}
  \label{eq:2}
  a_{\sym^2f}(m_1,m_2) \ = \
  \sum_{d|(m_1,m_2)}\lambda_{\sym^2f}\left(\frac{m_1}d,1\right)
  \lambda_{\sym^2f}\left(\frac{m_2}d,1\right)\mu(d),
\end{equation}
where $\mu$ is the M\"{o}bius function and
\begin{equation}
  \label{eq:3}
  \lambda_{\sym^2f}(r,1) \ = \  \sum_{s^2t=r}\lambda_f(t^2).
\end{equation}
The symmetric-square $L$-function of $f$ is then
\begin{equation}
  \label{eq:5}
  L(s,\sym^2f) \ = \  \sum_{m=1}^\infty \lambda_{\sym^2f}(m,1)m^{-s}.
\end{equation}
If, as before, $\alpha_p,\alpha_p^{-1}$ are the Satake parameters of~$f$, then
the parameters $\sigma_p{(j)}$ ($j=1,2,3$) of~$\sym^2f$ at any prime~$p$ are
the numbers $\alpha_p^2,1,\alpha_p^{-2}$.


Denoting by $\lambda_\phi(r)$ the $r$\textsuperscript{th} Hecke
eigenvalue of $\phi$, the Rankin-Selberg convolution
$L(s,\phi\times\sym^2f)$ is the Dirichlet series
\begin{align}
\label{eq:4}
L(s,\phi\times\sym^2f) & =
  \sum_{m_1,m_2\geq1}\lambda_\phi(m_1)\lambda_f(m_2)a_F(m_1,m_2)(m_1m_2^2)^{-s} \nonumber\\
  & =\  \sum_m \lambda_{\phi,\sym^2 f}(m)m^{-s},
\end{align}
 where
\begin{equation}
  \label{eq:13}
  \lambda_{\phi\times\sym^2f}(m) \ = \  \sum_{m_1m_2^2=m}\lambda_\phi(m_1)\lambda_f(m_2)a_F(m_1,m_2).
\end{equation}
In fact, also by the work of Kim-Shahidi~\cite{KiSh} (and the appendix by
Bushnell-Henniart), $L(s,\phi\times\sym^2f)$ is an automorphic $L$-function
$L(s,\pi)$ (for a suitable an automorphic representation $\pi$ of $\glsix$.)
This ensures the standard properties (entire of order one, bounded in vertical
strips, and functional equation) for $L(s,\phi\times\sym^2f)$.  In particular,
$L(s,\phi\times\sym^2f)$ conjecturally satisfies the Riemann Hypothesis in the
usual sense: $L(s,\phi\times\sym^2f)=0$ and $0\leq\Im s\leq 1$ implies $\Im
s=\frac12$.

The Satake parameters $\delta_p(j)$ ($j=1,\dots6$) of $\pi_p$ are
the six numbers $\alpha_p^{\pm2}\beta_p^{\pm1}$ and
$\beta_p^{\pm1}$.  Furthermore, each $L(s,\phi\times\sym^2f)$ has an
even functional equation.  The proof of this assertion is given in
Appendix~\ref{sec:appendix}.

For $\Re s$ large, the logarithmic derivative of
$L(s,\phi\times\sym^2f)$ is given by the Dirichlet series
\begin{equation}
  \label{eq:15}
  \frac{L'}L\left(s,\phi\times\sym^2f\right)
  \ = \ \sum_{m=0}^\infty \Lambda(m)a_{\phi\times\sym^2f}(m)m^{-s},
\end{equation}
where $\Lambda(m)$ is von Mangoldt's function and
\begin{equation}
  \label{eq:16}
  a_{\phi\times\sym^2f}(p^\nu) \ = \  \sum_{j=1}^6 \delta_p(j)^\nu.
\end{equation}

Define now the archimedean (gamma) factor
\begin{multline}
  \label{eq:6}
  L_\infty(s,\phi,\sym^2f)
  := \Gamma_{\R}(s+k-1+it_\phi)\Gamma_{\R}(s+k-1-it_\phi)\times \\
  \times \Gamma_{\R}(s+k+it_\phi)\Gamma_{\R}(s+k-it_\phi)
  \Gamma_{\R}(s+1+it_\phi)\Gamma_{\R}(s+1-it_\phi),
\end{multline}
where
\begin{equation}
  \label{eq:8}
  \Gamma_\R(s) := \pi^{-\frac s2}\Gamma\left(\frac s2\right).
\end{equation}
The completed $L$-function \be \Lambda(s,\phi\times\sym^2f) \ := \
L_\infty(s,\phi,\sym^2f)L(s,\phi\times\sym^2f)\ee  for
\eqref{eq:4} satisfies the functional equation
\begin{equation}
  \label{eq:10}
  \Lambda(s,\phi\times\sym^2f)
  \ = \  \Lambda(1-s,\phi\times\sym^2f).
\end{equation}
As the functional equation is even, we expect to observe either
$\soe$ or symplectic symmetry.

Following Rudnick and Sarnak \cite{RS}, we define the six
archimedean parameters $\mu_j$ ($j=1,\dots,6$) by the requirement
that $\frac12 + \mu_j$ is one of
\begin{equation}
  \label{eq:11}
  k\pm\frac12\pm it_\phi\quad\mbox{or}\quad \frac32\pm it_\phi.
\end{equation}

\subsubsection{Explicit Formula}
\label{sec:explicit-formula} A smooth form of the explicit formula
for $L(s,\phi\times\sym^2f)$ is as follows (see \cite{RS} for a
proof). Let $g\in C_c^\infty(\R)$ be an even Schwartz function
whose Fourier transform
\begin{equation}
  \hat g(y) \ = \  \int_{-\infty}^\infty g(x) e^{-2\pi i xy}\d x
\end{equation}
is compactly supported.  Let $R>0$ and write the non-trivial zeros
of $L(s,\phi\times\sym^2f)$ as $\rho_j=\frac12+i\gamma_j$; we have
$j\in\Z-\{0\}$ as the functional equation is even. Note $\gamma_j
\in \R$ is equivalent to GRH. Then
\begin{equation}
  \label{eq:14}
  \sum_j g\left(\frac{\gamma_j}{2\pi}\log R\right) \ = \
  \frac A{\log R} -
  2\sum_{p}\sum_{\nu=1}^\infty\hat g\left(\frac{\nu\log p}{\log R}\right)
  \frac{a_{\phi\times\sym^2f}(p^\nu)\log p}{p^{\nu/2}\log R},
\end{equation}
where
\begin{equation}
  \label{eq:1}
  A \ = \  \int_{-\infty}^{\infty} \sum_{j=1}^6 \left(\frac{\Gamma'_\R}{\Gamma_\R}
  \left(\mu_j+\frac12 + \frac{2\pi ix}{\log R}\right) + \frac{\Gamma'_\R}{\Gamma_\R}
  \left(\overline{\mu_j}+\frac12 + \frac{2\pi ix}{\log R}\right)\right)g(x)\d x.
\end{equation}

\subsubsection{Gamma Factor Contribution}\label{subsecgammasix}

Recall $\Gamma_\R(s) = \pi^{-s/2} \Gamma(s/2)$. Thus \be
\frac{\Gamma'_\R(s)}{\Gamma_\R(s)} \ = \  - \frac{\log \pi}{2} +
\foh \frac{ \Gamma'({\small\frac s2})}{\Gamma({\small\frac s2})}.
\ee Let $r=2\pi x/ \log R$. Then the sum in~\eqref{eq:1} equals
\be -6\log \pi + \foh \sum_{j=1}^6 \left[
\frac{\Gamma'}{\Gamma}\left(\fof + \frac{\mu_j}{2} +
\frac{ir}{2}\right) + \frac{\Gamma'}{\Gamma}\left(\fof +
\frac{\overline{\mu_j}}{2} - \frac{ir}{2}\right) \right], \ee
where $\mu_j = k \pm \foh \pm it_\phi$ (for four values) and
$\frac{3}{2} \pm it_\phi$ (for the other two). We use (see
\cite{ILS} or \cite{GR} 8.363.3) that for $a,b\in\R$, $a>0$, \be
\frac{\Gamma'}{\Gamma}\left(a + bi\right) +
\frac{\Gamma'}{\Gamma}\left(a - bi\right) \ = \
2\frac{\Gamma'}{\Gamma}\left( a \right) + O(a^{-2} b^{2}), \ee and
for $\alpha \geq \fof$, \be \frac{\Gamma'}{\Gamma}\left(\alpha +
\fof\right)  \ = \  \log \alpha + O(1). \ee Thus, in the
$\Gamma'/\Gamma$ factors, the $\mu_j = \frac{3}{2} \pm it_\phi$
terms are $O(1)$ with respect to $k$. Set $a_+ = \foh$ and $a_- =
0$. Matched in complex-conjugate pairs, the other eight terms give
\bea & & \frac{\Gamma'}{\Gamma}\left(\frac{k}{2} + a_\pm
+i\left(\pm
      \frac{t_\phi}{2}+r\right) \right) +
  \frac{\Gamma'}{\Gamma}\left(\frac{k}{2} + a_\pm -i\left(\pm
      \frac{t_\phi}{2}+r\right) \right) \nonumber\\
  & & \ \ \ \ \ \ \ \ \ \ \ = \
  2\frac{\Gamma'}{\Gamma}\left( \frac{k}{2} + a_\pm\right) + O\left(
    \frac{|t_\phi|^2 + r^2}{k^2} \right) \eea and
\bea\label{eq:21}
    \frac{\Gamma'}{\Gamma}\left( \frac{k}{2} + a_\pm\right) \ = \
    \log\left( \frac{k}{2} + a_\pm - \fof\right) + O(1)  \ = \  \log k + O(1).
\eea Note for $k \geq 2$, the condition of having the argument
greater than $\fof$ is trivially met.  The main term in the sum
in~\eqref{eq:1} is simply $\frac12\cdot4\cdot2\log k=4\log k$. The
main contribution of the term $\frac{A}{\log R}$ in~\eqref{eq:14}
is
\begin{equation}
  \label{eq:22}
  \frac{4\log k }{\log R}\int_{-\infty}^\infty g(x)dx
   \ =\
    \frac{\log k^4}{\log R}
  \cdot \hat g(0).
\end{equation}
In the $k$-aspect, $\phi\times\sym^2 f$ looks like a $\glfour$
object. A natural choice for the analytic conductors is therefore
$k^4$. With this scaling of the zeros, the test function on the
the left-hand side of~\eqref{eq:22} is evaluated at points which
have mean average spacing one near the central point (Riemann's
classical critical zero-counting formula).  As the quotient
depends only on the logarithm of the conductor, as $k \to \infty$
the choice of any fixed constant multiple of $k^4$ for the
conductor will give the same answer (see \cite{ILS}). We have
proved

\begin{lem}\label{lemcontributiongammafactors} For $L(s,\phi\times\sym^2
f)$, up to lower order terms the contribution from the
$\Gamma$-factors in the explicit formula equals $\widehat{g}(0)$,
and the analytic conductor equals $k^4$. Assuming GRH, the
non-trivial zeros of $L(s,\phi\times\sym^2 f)$ are $\foh +
i\gamma_{\phi\times\sym^2 f}^{(j)}$ with $\gamma_{\phi\times\sym^2
f}^{(j)} \in \R$. Taking $R = k^4$, the explicit
formula~\eqref{eq:14} becomes
\begin{multline}
\label{eqexplicitformulafinal}
\sum_j g\left(\frac{\gamma_{\phi\times\sym^2 f}^{(j)}}{2\pi}\log R\right) \ = \  \\
\hg{0} - 2\sum_{p}\sum_{\nu=1}^\infty\hat g\left(\frac{\nu\log
p}{\log R}\right) \frac{a_{\phi\times\sym^2f}(p^\nu)\log
p}{p^{\nu/2}\log R} + O\left(\frac1{\log
    R}\right).
\end{multline}
\end{lem}
The appearance of the term $\hg{0}=\int g(x)\d x$ on the right-hand
side of~\eqref{eqexplicitformulafinal} naturally corresponds to the
(expected) term $\int \hg{\xi}\delta(\xi)\d\xi$ due to the delta
mass at the origin in the Fourier transform of the $1$-level density
(see \eqref{eq:FTphi1ld}).  The second term (double sum) above will
eventually be matched to $\int\hg{\xi}\frac12\eta(\xi)\d\xi$, and
this will exclude symplectic as a possibility.

\begin{rek} It is fortunate for us that the analytic conductors of
$L(s,\phi \times \sym^2 f)$ depend weakly on $f$. Specifically, as
the only dependence on $f$ is through its weight $k$, one scaling
works for all elements of our family. Oscillating conductors in a
family can sometimes be handled (one recourse is to use the average
log conductor as in \cite{Si,Yo}; another approach is a more careful
analysis and sieving, as in \cite{Mil2} where the conductors are
monotone).
\end{rek}

\subsubsection{Relation of $a_{\phi\times\sym^2f}$ to $\lambda_f$ and
  $\lambda_\phi$}
\label{sec:a_phi-expansion} To evaluate the double sum
in~\eqref{eqexplicitformulafinal}, we express
$a_{\phi\times\sym^2f}(p^\nu)$ in terms of
$\lambda_{f},\lambda_\phi$. Note
\begin{eqnarray}\label{eq:24}
    \apsf(p^\nu)
    &\ = \ & \alpha_p^{2\nu}\beta_p^\nu + \alpha_p^{-2\nu}\beta_p^\nu + \alpha_p^{2\nu}\beta_p^{-\nu} +
    \alpha_p^{-2\nu}\beta_p^{-\nu}+ \beta_p^\nu+\beta_p^{-\nu} \nonumber\\
    &\ = \ & (\alpha_p^{2\nu} + 1 + \alpha_p^{-2\nu})(\beta_p^\nu + \beta_p^{-\nu}).
\end{eqnarray}

\bigskip

\noindent\textbf{Case $\nu=1$:} \label{sec:a_nu_1} We  have
\begin{equation}
  \label{eq:17}
  \begin{split}
    \apsf(p) &\ = \ (\alpha_p^2+1+\alpha_p^{-2})(\beta_p+\beta_p^{-1})
    \ = \ \lambda_f(p^2)\lambda_\phi(p).
  \end{split}
\end{equation}

\bigskip

\noindent \textbf{Case $\nu\geq2$:} \label{sec:a_nu_geq_2} We
have
\begin{multline}
  \label{eq:26}
  \alpha_p^{2\nu} + 1 + \alpha_p^{-2\nu} \ = \  (\alpha_p^{2\nu} + \alpha_p^{2(\nu-1)} + \dots +
  \alpha_p^{-2(\nu-1)} + \alpha_p^{-2\nu}) \\
  - (\alpha_p^{2(\nu-1)}+\dots+\alpha_p^{-2(\nu-1)}) + 1,
\end{multline}
and
\begin{multline}
  \label{eq:27}
  \beta_p^\nu + \beta_p^{-\nu}
  \ = \  \beta_p^\nu + \beta_p^{\nu-2} + \dots + \beta_p^{-(\nu-2)} + \beta_p^{-\nu} \\
   - (\beta_p^{\nu-2} + \dots + \beta_p^{-(\nu-2)}),
\end{multline}
yielding
\begin{equation}
  \label{eq:25}
  \begin{split}
    \apsf(p^\nu) &\ = \  (\alpha_p^{2\nu} + 1 + \alpha_p^{-2\nu})(\beta_p^\nu + \beta_p^{-\nu}) \\
    &\ = \ (\lambda_f(p^{2\nu}) - \lambda_f(p^{2(\nu-1)})+1)(\lambda_\phi(p^\nu)-\lambda_\phi(p^{\nu-2})).
  \end{split}
\end{equation}

\bigskip

Of course $\lambda_\phi(p^{\nu-2})=1$ when $\nu=2$.

\subsubsection{Summary}

We have shown

\begin{lem}\label{lemexpansionafpafpp}
  \begin{align}
    \label{eq:29}
    \afsp(p) &\ = \  \glp(p) \lambda_f(p^2) \\
    \label{eq:30}
    \afsp(p^2) &\ = \  (\glp(p^2) - 1) \cdot (\lambda_f(p^4)-\lambda_f(p^2)+1).
  \end{align}
\end{lem}
As we shall see below, the single term `$-1$' in the first factor of
\eqref{eq:30} is responsible for flipping  the symmetry from
symplectic (for the $\{\sym^2f\}$ family of \cite{ILS} which had
$a_{\sym^2 f} = \lambda_f(p^4)-\lambda_f(p^2)+1$) to $\soe$ (for the
$\{\phi\times\sym^2f\}$ family we are considering). This behavior is
described in more detail in \cite{DM}.

Using the results from Kim-Sarnak~\cite{K}, we have $|\beta_p^{\pm
1}| \leq p^{\frac{7}{64}}$. Since $|\alpha_p| \leq 1$, equation
\eqref{eq:16} yields
\begin{equation}
\label{eq:32}
  \afsp(p^\nu)\ =\  (\beta_p^\nu + \beta_p^{-\nu}) \cdot
  (\alpha_p^{2\nu} + \alpha_p^{-2\nu} + 1)\ \ll\ p^{\frac{7\nu}{64}}.
\end{equation}
Therefore \be \frac{\afsp(p^\nu)}{p^{\nu/2}}\ \ll\ p^{-
\frac{25\nu}{64}}. \ee

This immediately implies

\begin{lem}\label{lemnocontributionnu} The contribution from terms with $\nu\geq3$ in
\eqref{eqexplicitformulafinal} can be absorbed into the error
term.
\end{lem}

\begin{rek} We do not need the full strength of $|\beta_p^{\pm 1}| \ll
p^{\frac7{64}}$; any exponent less than $\frac{1}{6}$ suffices.
Without such a bound, we would later need to obtain cancellation
when averaging the Fourier coefficients over the family (a result
of this nature is significantly weaker than proving bounds towards
Ramanujan, and follows from the Petersson formula).
\end{rek}


\subsection{$1$-Level Density}\label{subseconelevelsix}

As $\FD_{\phi\times \sym^2 H_k} = \{\phi \times \sym^2 f, f \in
H_k\}$, we have $|\FD_{\phi\times \sym^2 H_k}| = |H_k|$. For each
$L$-function from $\FD_{\phi\times \sym^2 H_k}$ we calculate the
$1$-level density for its low lying zeros via the explicit
formula; we then average over the family $\FD_{\phi\times \sym^2
H_k}$.

\subsubsection{Preliminaries}\label{sec:prelimsgl6case}

Let $g$ be an even Schwartz function with ${\rm supp}(\widehat{g})
\subset (-\sigma,\sigma)$.  Following Iwaniec-Luo-Sarnak \cite{ILS}
or Royer \cite{Ro}, we consider a weighted average over the
family~$\FD_{\phi\times \sym^2 H_k}$ of the
expressions~\eqref{eqexplicitformulafinal}.  The weight factors that
we use are $\omega_f = \zeta(2)/L(1,\phi\times\sym^2f)$. These are
positive (by GRH for $L(s,\sym^2f)$), slowly varying, and satisfy
\begin{equation}
  \label{eq:52}
  k^{-\epsilon}\ \ll_\epsilon\ \frac1{L(1,\sym^2f)}\ \ll_\epsilon\ k^\epsilon
\end{equation}
and
\begin{equation}
  \label{eq:53}
  \frac{1}{|\hkn|}\sum_{f\in H_k}\frac{\zeta(2)}{L(1,\sym^2f)}\ =\ 1+O\left(\frac1k\right).
\end{equation}
To simplify the application of the Petersson Formula, we have
introduced the slowly varying weights $\zeta(2)/L(1,\sym^2 f)$;
arguing along the lines of \cite{ILS} allows one to remove these
weights at no cost. We have chosen to leave in the weights in order
to emphasize the features of this $\glsix$ family.

By GRH, we may denote the non-trivial zeros of $L(s,\phi\times
\sym^2 f)$ by $\foh + i\gamma_{\phi\times\sym^2 f}^{(j)}$ with
$\gamma_{\phi\times\sym^2 f}^{(j)}\in\R$. Let $R=k^4$. All
$L(s,\phi\times\sym^2 f)$ have the same analytic conductor, which
up to lower order terms is $k^4$. Averaging
\eqref{eqexplicitformulafinal} by incorporating the weights and
using Lemma \ref{lemnocontributionnu} to absorb the $\nu \ge 3$
terms into the error shows that the $1$-level density for the
family $\FD_{\phi\times \sym^2 H_k}$ is \bea\label{eq:33} & &
D_{1,\FD_{\phi\times \sym^2 H_k}}(g) \nonumber\\ & = &
\frac{1}{|\hkn|} \sum_{f \in
  \hkn}\frac{\zeta(2)}{L(1,\sym^2 f)} \sum_j g\left(
  \gamma_{\phi\times \sym^2 f}^{(j)} \frac{\log c_{f_\phi}}{2\pi}\right) \nonumber\\
& = &  \hg{0} - \frac{2}{|\hkn|}
  \sum_{f \in \hkn}\frac{1}{L(1,\sym^2f)} \sum_{\nu=1}^2 \sum_{p=2}^{R^\sigma}\frac{\afsp(p^\nu)\log
p}{p^{\nu/2}\log R} \widehat{g}\left(\nu \frac{\log p}{\log
R}\right)\nonumber\\ & & \ \ \ \ +\ O\left(\frac1{\log
    R}\right).
\eea

We are left with analyzing the contribution from the $\nu=1,2$
terms, and comparing this to \eqref{eq:FTphi1ld}. For small
support, we will show there is no contribution from the $\nu = 1$
term, and the $\nu = 2$ term contributes $\foh g(0)$. This proves
the symmetry group is neither unitary nor symplectic. We cannot
discard the $\so$ and $\soo$ symmetries; however, we will be able
to eliminate them later by studying the $2$-level density.

\begin{rek}
  Since Iwaniec-Luo-Sarnak exclusively use $1$-level density
  arguments, they must use extra averaging to extend their support
  past $[-1,1]$.  By studying the $2$-level density we provide
  compelling evidence for the underlying symmetry being
  $\soe$ without extra averaging.
\end{rek}




\subsubsection{Contribution from $\nu = 1$}\ \\

We must evaluate \be T_1 \ = \ \frac{1}{|\hkn|} \sum_{f\in
  \hkn}\frac{\zeta(2)}{L(1,\sym^2 f)}\sum_{p=2}^{R^\sigma}
\frac{\afsp(p)\log p}{\sqrt{p}\log R}\ \widehat{g}\left(\frac{\log
    p}{\log R}\right). \ee By Lemma \ref{lemexpansionafpafpp},
$\afsp(p) = \glp(p) \lambda_f(p^2)$. Since $\lambda_{f}(1)=1$,
\be T_1 \ = \ \sum_{p=2}^{R^\sigma} \frac{\glp(p)\log
p}{\sqrt{p}\log R} \widehat{g}\left(\frac{\log
    p}{\log R}\right) \frac{\zeta(2)}{|\hkn|} \sum_{f\in
  \hkn}\frac{\lambda_f(1)\lambda_f(p^2)}{L(1,\sym^2
  f)}. \ee
We are led to studying \be \frac{\zeta(2)}{|\hkn|} \sum_{f\in
\hkn} \frac{\lambda_f(1)\lambda_{f}(p^2)}{L(1,\sym^2 f)}. \ee

We have a non-diagonal term since $p^2\neq1$. By
\eqref{eqilstwothree} of Lemma \ref{lempeterssonfromils},
$\delta(1,p^2) = 0$; for $p \ll k$ these terms are $\ll
\frac{p}{2^k}$. Substituting into the expansion for $T_1$, we see
there is no contribution for $R^{\sigma} < k$. Since $R = k^4$, this
implies there is no contribution for $\sigma < \frac{1}{4}$. Note
that in executing the prime sum, any polynomial bound on $\glp(p)$
suffices, since the decay in $k$ is exponential.

For primes $p > R^{\fof}$, we cannot use \eqref{eqilstwothree};
instead we use \eqref{eqilstwotwo}, which gives $\ll
\frac{\sqrt{p}\log p}{k^{5/6}}$. This yields a $p$-sum of \be
\frac{1}{k^{5/6}} \sum_{p}^{k^{4\sigma}} \frac{\log^2 p}{\log R}
\frac{\glp(p)\sqrt{p}}{\sqrt{p}}. \ee Let $\delta$ be the the best
bound towards Ramanujan for $\glp(p)$; namely, $\glp(p) \ll
p^\delta$ (the Ramanujan conjecture is $\delta = 0$). We find this
sum is $\ll k^{4\sigma  (1+\delta) - \frac{5}{6}}$; thus, $\sigma <
\frac{5}{24   (1 + \delta)}$. Even assuming Ramanujan does not help
---this bound is worse than the previous one. Thus
\eqref{eqilstwothree} is better, and we obtain that there is no
contribution for support up to $\fof$.

\begin{rek} The reason there is no contribution for small
support is that we have a non-diagonal term in the the Petersson
formula. \end{rek}

\subsubsection{$\nu = 2$}\label{subsec:1ldnu2}

We must evaluate \be T_2  \ = \  -\frac{2}{|\hkn|} \sum_{f\in
\hkn}\frac{\zeta(2)}{L(1,\sym^2 f)}\sum_{p=2}^{R^\sigma}
\frac{\afsp(p^2)\log p}{p\log R}\ \widehat{g}\left(2\frac{\log
p}{\log R}\right). \ee By Lemma \ref{lemexpansionafpafpp},
$\afsp(p^2) = (\glp(p^2) - 1) \cdot
(\lambda_{f}(p^4)-\lambda_f(p^2)+1)$. Almost all of the terms are
non-diagonal. Using $\lambda_f(1) = 1$ we have the following terms:
from $\glp(p^2)$, we get \bea & & \glp(p^2)\cdot
\lambda_f(1)\lambda_f(p^4), \ \ \
-\glp(p^2)\cdot\lambda_f(1)\lambda_f(p^2), \ \ \
\glp(p^2)\cdot\lambda_f(1)\lambda_{f}(1). \nonumber\\ \eea

The first two terms are non-diagonal; the Petersson formula yields
no contribution for small support. The third term \emph{is}
diagonal. For small support, up to lower order terms it yields $+1$
by \eqref{eq:53}: \be \frac{\zeta(2)}{|\hkn|}\sum_{f \in \hkn}
\frac{\lambda_f(1)\lambda_f(1)}{L(1,\sym^2 f)} \ = \ 1 +
O\left(\frac1k\right). \ee This gives \be\label{eqglsixsymphi}
-2\sum_{p=2}^{R^\sigma} \frac{\glp(p^2)\log p}{p\log R}\
\widehat{g}\left(2\frac{\log p}{\log R}\right). \ee By GRH for
$L(s,\sym^2 \phi)$, this sum is $O(\frac{1}{\log R})$ (see Section 4
of \cite{ILS}).

We now handle the three terms from the $-1$ in the first factor of
$\afsp(p^2)$; these are \be - \lambda_f(1)\lambda_{f}(p^4), \ \ \
\lambda_f(1)\lambda_f(p^2),\ \ \ -\lambda_f(1)\lambda_f(1). \ee The
first two are non-diagonal, and by the Petersson formula do not
contribute for small support. The third term, however, \emph{is} a
diagonal term; up to lower order corrections, from the Petersson
formula its contribution is $-1$, and we are left with \be 2\sum_p
\frac{\log p}{p\log R}\ \hg{2\frac{\log p}{\log R}}. \ee By Lemma
\ref{thmprimesums}, the above sum (up to lower order terms) is $\foh
g(0)$.

Therefore, for small support, the $\nu = 2$ piece contributes $\foh
g(0) + o(1)$, with the main term arising from the sixth term in the
expansion of $\afsp(p^2)$.  At this point we have enough evidence to
discard the unitary and symplectic symmetries. Since the functional
equations in \eqref{eq:10} are even, this certainly points to the
underlying symmetry being $\soe$, but we cannot yet discard the full
orthogonal or $\soo$ symmetries. This will be done in
\S\ref{sec:2-level-density}.

We now determine how large we may take the support. Of the six
pieces which do not contribute to the main term, the worst error
term is from $\glp(p^2) \cdot \lambda_{f}(1)\lambda_f(p^4)$. By
\eqref{eqilstwothree}, if $1\cdot p^4 \ll k^2$, the sum over $f\in
\hkn$ is $\ll \frac{p^2}{2^k}$. Again, any polynomial bound for
$\glp(p^2)$ yields the sum over primes $p \ll k^\foh$ is a lower
order term. Since $R = k^4$, this yields no contribution for
$\sigma < \frac{1}{8}$. Therefore, up to lower order terms the
contribution is $\foh g(0)$ for $\sigma < \frac18$.

The reason for the sharp decrease in support (relative to the $\nu =
1$ term) is because we have a $\lambda_f(1)\lambda_{f}(p^4)$.
Another possibility is to use \eqref{eqilstwotwo}, which gives for
the $\glp(p^2) \cdot \lambda_f(1)\lambda_f(p^4)$ term
\be\label{eqsecondtry} \frac{1}{k^{5/6}} \sum_{p}^{k^{4\sigma}}
\frac{\log^2 p}{\log R} \frac{\glp(p^2)\cdot p}{p}. \ee From
\eqref{eqlambdaexpand}, $\gl_\phi(p^2) = \gl_\phi(p)^2 - 1$.
Substituting into \eqref{eqsecondtry}, the -1 does not contribute
for $\sigma < \frac5{24}$, and we are left with bounding
\be\label{eqlambdaexpandtwo} \frac{1}{k^{5/6}}
\sum_{p=2}^{k^{4\sigma}} \frac{\log^2 p}{\log R}\ \glp(p)^2 \ \ll \
\frac{\log k^4}{k^{5/6}} \sum_{n=1}^{k^{4\sigma}} |\glp(n)|^2. \ee
One could use bounds towards Ramanujan; however, all we need is that
Ramanujan holds on average, namely the sum of $|\gl_\phi(n)|^2$ (see
\cite{Iw1}, equation 8.7) is \be\label{eqIwX} \sum_{n=1}^X
|\gl_\phi(n)|^2 \ \ll_{\phi} \ X. \ee This yields the sum in
\eqref{eqlambdaexpandtwo} does not contribute for $\sigma <
\frac5{24}$.

A similar argument applied to the other terms show none of them
contribute as well for such support; one must check that the error
term for the $\gl_f(1)\gl_f(1)$ term (which is the diagonal piece
responsible for $\foh g(0)$) does not contribute in this range.
This completes the proof of the first part of Theorem
\ref{thm:soefam}.

As  $\text{supp}(\hphi) \subset (-1,1)$, while every
$\Lambda(s,\phi\times\sym^2 f)$ has even functional equation, we
cannot conclude the symmetry is $\soe$ and not $\so$ or $\soo$.
There are two natural ways to try and increase the support. The
first is to average over even Maass forms $\phi$; unfortunately,
we would have to let $t_j$ grow to a power of $k$, which would
change the conductor arguments. Another approach is to average
over the weight $k$ (as in the investigation of $\sym^2f$ in
\cite{ILS}). By averaging over weight, they triple the support;
however, as we start with support less than $\fot$, such methods
are insufficient to break $(-1,1)$. We therefore study the 2-level
density, which even for arbitrarily small support  can distinguish
the three orthogonal groups (see \cite{Mil1,Mil2}).


\subsection{$2$-Level Density}\label{sec:2-level-density}

We complete the proof of Theorem \ref{thm:soefam}. As in
\S\ref{subseconelevelsix}, for convenience in applying the
Petersson formula we study a weighted 2-level density. Thus
\eqref{eq:7b} becomes
\begin{eqnarray}\label{eq:28}
& &  D_{2,\mathcal{F}_{\phi\times\sym^2 H_k}}(g)  \ = \ \nonumber\\
& & \frac1{|H_k|} \sum_{f\in H_k}
  \frac{\zeta(2)}{L(1,\sym^2 f)} \sum_{j_1,  j_2 \atop j_1 \neq
    \pm j_2} g_1\left(\frac{\log R}{2\pi}\gamma_f^{(j_1)}\right)
  g_2\left(\frac{\log R}{2\pi}\gamma_f^{(j_2)}\right).\ \ \ \ \ \
  \ \
\end{eqnarray}
The sum is over zeros $j_1 \neq \pm j_2$. As all functional
equations are even and we are assuming GRH, the zeros occur in
complex conjugate pairs; this is the first time the sign of the
functional equation enters our arguments. We may rewrite the
$2$-level expression as a sum over all pairs of zeros, minus twice
the sum over all zeros, yielding
\begin{eqnarray}
\label{eq:31} & & D_{2,\mathcal{F}_{\phi\times \sym^2 H_k}}(g) \ =
\ \nonumber\\ & & \frac1{|H_k|} \sum_{f\in H_k}
\frac{\zeta(2)}{L(1,\sym^2f)} \sum_{j_1,j_2}g_1\left(\frac{\log
R}{2\pi}\gamma_f^{(j_1)}\right) g_2\left(\frac{\log
R}{2\pi}\gamma_f^{(j_2)}\right) \nonumber\\ & & \ -
2D_{1,\mathcal{F}_{\phi\times \sym^2 H_k}}(g_1g_2).
\end{eqnarray}
In the above, the second term is a $1$-level density with test
function $(g_1g_2)(x)$. For small support, we have shown this term
is just $\widehat{g_1g_2}(0) + \foh (g_1g_2)(0)$. Crucial in the
above expansion is that each $\Lambda(s,\phi\times\sym^2 f)$ has
even sign. This allows us to pair off the zeros, $\gamma_f^{(j)}$
and $\gamma_f^{(-j)}$. Thus summing over distinct zeros is the same
as subtracting off twice a 1-level sum over all zeros. This would be
false if the functional equation were odd. In that case we would
have to add back $g_1(0)g_2(0)$ for the extra zero at the central
point, and in fact the presence or absence of this additional term
is the cause of the differences in the 2-level densities of the
three orthogonal groups.

For the first term in \eqref{eq:31}, since we are summing over all
zeros, we may use the explicit formula for the sum over each
$j_i$. Let
\begin{eqnarray}
\label{eq:38} \bfsp(p) &  \ = \ & \afsp(p) \nonumber\\
\bfsp(p^2) & \ = \ & \afsp(p^2) + 1.
\end{eqnarray}
In the expansions with the explicit formula, we isolate the
contribution from the part of the $\nu = 2$ term which contributes
$\foh g_i(0)$ for small support; this is the $+1$ term in
$\bfsp(p^2)$. We have also removed the error terms arising from
$\nu \geq 3$; the argument is standard (see \cite{Ru2,RS}). We are
left with considering the weighted average of
\begin{equation}
\label{eq:37} \prod_{i=1}^2 \left[ \left(\hgi{0} + \foh
g_i(0)\right)- 2\sum_{\nu_i=1}^2\sum_{p_i}
\frac{\bfsp(p^{\nu_i})\log p}{p^{\nu_i/2}\log R} \hgi{\nu_i
\frac{\log p_i}{\log R}}\right].
\end{equation}

There are three terms for each $i$. For small support, we have
shown the $\nu_i$-sums by themselves do not contribute, though we
will see that there are contributions when a $\nu_1$-sum hits a
$\nu_2$-sum. Thus when a $\hgi{0} + \foh g_i(0)$ hits a $\nu$-sum,
there is no contribution. We have $[\hgo{0}+\foh g_1(0)]\cdot
[\hgt{0}+\foh g_2(0)]$, plus the weighted average of the mixed
sums \be 4 \sum_{p_1} \sum_{p_2} \frac{ \bfsp(p_1^{\nu_1})
\bfsp(p_2^{\nu_2}) \log p_1 \log p_2}{p_1^{\nu_1/2} p_2^{\nu_2/2}
\log^2 R} \hgo{\nu_1 \frac{\log p_1}{\log R}} \hgt{\nu_2
\frac{\log p_2}{\log R}} \ee for $(\nu_1,\nu_2) \in \{(1,1),
(2,1), (2,1), (2,2)\}$. For each pair there are two cases, when
$p_1 = p_2$ and $p_1 \neq p_2$.  Since we are only interested in
the $2$-level density for arbitrarily small support as a means to
distinguish $\soe$ from orthogonal and $\soo$ symmetry, we do not
record how large we may take the support.

\subsubsection{$(1,1)$ Terms}
\label{sec:1-1-terms} We have \bea & & \frac4{|H_k|} \sum_{f \in
\hkn} \weight \sum_{p_1} \sum_{p_2} \frac{ \afsp(p_1)
\afsp(p_2) \log p_1 \log p_2}{\sqrt{p_1 p_2} \log^2 R} \nonumber\\
& & \times \hgo{\frac{\log p_1}{\log R}} \hgt{\frac{\log p_2}{\log
R}}. \eea As $\afsp(p) = \glp(p) \lambda_f(p^2)$, we have
\begin{multline}
 \frac4{|H_k|} \sum_{f \in \hkn} \weight \sum_{p_1} \sum_{p_2}
\glp(p_1)\glp(p_2)\frac{ \lambda_f(p_1^2) \lambda_{f}(p_2^2) \log
p_1 \log
p_2}{\sqrt{p_1 p_2} \log^2 R} \\
\times \hgo{\frac{\log p_1}{\log R}} \hgt{ \frac{\log p_2}{\log
R}}.
\end{multline}
If $p_1 \neq p_2$, when we use the Petersson formula there is no
contribution for small support, since it is a non-diagonal term.
Note $p_1 = p_2$ is a diagonal term, giving $\frac1{|H_k|} \sum_f
\weight\lambda_f(p^2) \lambda_f(p^2) = 1 + O(\sqrt{p^4}/2^k)$. For
sufficiently small support the error term is negligible, and thus
the diagonal term is \be 4\sum_p \glp(p) \glp(p) \frac{\log^2
p}{p\log^2 R} \hgo{\frac{\log p_1}{\log R}}\ \hgt{
  \frac{\log p_2}{\log R}} \ + \ o(1).
\ee We use $\glp(p)\glp(p) = 1 + \glp(p^2)$; we saw in
\S\ref{subsec:1ldnu2} that the $\glp(p^2)$ term will not
contribute by GRH for $L(s,\sym^2 \phi)$. The $+1$ will
contribute, with test function $\widehat{g_1}\widehat{g_2}$. By
Lemma \ref{thmprimesums}, we have \be 4\sum_p \frac{\log^2
p}{p\log^2 R} \cdot \widehat{g_1}\widehat{g_2}\left( \frac{\log
p}{\log R}\right)\ =\
 2 \int |u| \widehat{g_1}\widehat{g_2}(u)\d u + O\left(\frac1{\log R}\right). \ee

\subsubsection{$(1,2)$ and $(2,1)$ Terms}
\label{sec:1-2-2-1-terms} As these terms are handled identically,
we confine ourselves to $(1,2)$. We have weighted averages of \be
 4 \sum_{p_1} \sum_{p_2} \frac{ \afsp(p_1)\bfsp(p_2^2)\log p_1
\log p_2}{p_1^{1/2} p_2 \log^2 R} \hgo{\frac{\log p_1}{\log R}}
\hgt{2 \frac{\log p_2}{\log R}}. \ee

By Lemma \ref{lemexpansionafpafpp} and \eqref{eq:38}
\begin{align}
 \afsp(p) & \ = \  \glp(p) \lambda_f(p^2) \\
 \bfsp(p^2) & \ = \  \glp(p^2)(\lambda_f(p^4)-\lambda_f(p^2)+1) - (\lambda_f(p^4)-\lambda_f(p^2)).
\end{align}
If $p_1 \neq p_2$, all terms are non-diagonal, and the Petersson
formula yields no contribution for small support.

If $p_1 = p_2$, only two terms give diagonal terms:
$\glp(p)\glp(p^2) \cdot \lambda_f(p^2)\lambda_{f}(p^2)$ and
$-\glp(p) \lambda_f(p^2)\lambda_{f}(p^2)$. However, while the
Petersson formula will give a $\pm 1$ for each of these diagonal
terms, this is immaterial since we are dividing by $p^{\frac32}$.
Using the Kim-Sarnak bound of $p^{\frac{7}{64}}$ for Maass forms
is enough to show there is no contribution, for any support, from
these terms. We are dividing by $p^{\frac{3}{2}}$, and we have at
most $p^{\frac{7\cdot 3}{64}}\log p$ in the numerator. This gives
$p^{-1-\frac{11}{64}} \log p$, which yields a $O(\frac{1}{\log
R})$ contribution. Arguing as in \eqref{eqsecondtry} to
\eqref{eqIwX}, one may replace the Kim-Sarnak bound with any
non-trivial bound towards Ramanujan.

\subsubsection{$(2,2)$ Term}
\label{sec:2-2-term} We now consider the $(2,2)$ term. We have a
weighted average of  \be 4 \sum_{p_1} \sum_{p_2} \frac{
\bfsp(p_1^2) \bfsp(p_2^2) \log p_1 \log p_2}{p_1 p_2 \log^2 R}
\hgo{2 \frac{\log p_1}{\log R}} \hgt{2 \frac{\log p_2}{\log R}},
\ee where by Lemma \ref{lemexpansionafpafpp} and \eqref{eq:38}
\bea \bfsp(p^2) & \ =
\ & \glp(p^2) \cdot (\lambda_f(p^4)-\lambda_{f}(p^2)+1)\nonumber\\
& & - 1 \cdot(\lambda_{f}(p^4)-\lambda_f(p^2)). \eea When $p_1 =
p_2$, as $\glp(p) = O(p^{\delta})$ for some $\delta \in [0,\foh]$,
$\frac{\bfsp(p^2)^2}{p^2} = O(p^{4\delta-2})$, and these terms
will not contribute for any $\delta < \fof$. We do not need the
full strength of the Kim-Sarnak bound; the $\frac5{28}$ of
\cite{BDHI} suffices.

We are left with the case $p_1 \neq p_2$. The only diagonal term
will be $\glp(p_1^2) \glp(p_2^2) \cdot
\lambda_{f}(1)\lambda_f(1)$. For small support, the other terms
will not contribute, and by Petersson's formula we  have \be
4\prod_{i=1}^2 \sum_{p_i} \frac{\glp(p_i^2) \log p_i}{p_i \log R}
\hgi{2 \frac{\log p_i}{\log R}}. \ee As in \S\ref{subsec:1ldnu2},
by GRH for $L(s,\sym^2 \phi)$ each prime sum is $O(\frac{1}{\log
R})$. Thus there is no contribution from the $(2,2)$ terms.

\subsubsection{Summary} \label{sec:summary} We have shown
\begin{multline}
  \label{eq:34}
  D_{2,\mathcal{F}_{\phi\times \sym^2 H_k}}(g) \ = \  \left[\hgo{0}+\foh g_1(0)\right]
  \cdot\left[\hgt{0}+\foh g_2(0)\right] \\
  + 2\int_u |u|\widehat{g_1}\widehat{g_2}(u)\d u - 2\widehat{g_1g_2}(0)
  - g_1(0)g_2(0),
\end{multline}
and \eqref{eq:34} agrees only with the $2$-level density for
$\soe$ (see \eqref{eq:soetwo}), completing the proof of Theorem
\ref{thm:soefam}.


\section{$\mathcal{F}_{\phi\times H_k} = \{\phi\times f: f \in H_k\}$}\label{sec:L-phi-f}

In this section we prove Theorem \ref{thm:sympfam}, namely that the
symmetry type of the family $\mathcal{F}_{\phi\times H_k} =
\{\phi\times f: f \in H_k\}$ ($k\to\infty$) agrees only with
symplectic, where $\phi$ is a fixed even cuspidal Hecke-Maass form
with eigenvalue $\lambda_\phi=\frac14+t_\phi^2$ and $f\in\hkn$ is a
Hecke holomorphic modular form of weight~$k$. These $L$-functions
have associated Euler products of degree~4 and are indeed associated
to automorphic representations of $\glfour$~\cite{Ra}.  As in
\S\ref{sec:defn-gamma-fe} we first derive the explicit formula and
find the analytic conductors for the family. In order to compute the
$1$-level density we analyze the local parameters and find evidence
for symplectic symmetry (for small support). Since the symplectic
$1$-level density is distinguishable from the other classical
compact groups for arbitrarily small support, there is no need to
investigate the $2$-level density. As the arguments are similar to
those for $\phi\times\sym^2f$, we merely sketch the calculations
below. See also Appendix \ref{sec:appendix} for details of the
determination of the gamma factors and signs of the functional
equation.

\subsection{Logarithmic Derivative, Gamma Factors, Functional
  Equation} \label{sec:Log-gamma-FE} In terms of the Fourier coefficients
$\{\lambda_f(n)\},\{\lambda_\phi(n)\}$ and the Maass eigenvalue
$t_\phi$ ($\frac14+t_\phi^2$ is the Laplacian eigenvalue), we have
\begin{equation}
  \label{eq:39}
  \begin{split}
    L(s,\phi\times f) & = \zeta(2s)\sum_m \lambda_\phi(m)\lambda_f(m)m^{-s} =
    \sum_m \lambda_{\phi\times f}(m) m^{-s},
  \end{split}
\end{equation}
where
\begin{equation}
  \label{eq:40}
  \lambda_{\phi\times f}(m) \ = \ \sum_{m_1^2 m_2=m}\lambda_\phi(m_2)\lambda_f(m_2).
\end{equation}
The logarithmic derivative of $L(s,\phi\times f)$ is
\begin{equation}
  \label{eq:45}
  \frac{L'}L\left(s,\phi\times f\right)\ = \ \sum_m \Lambda(m)a_{\phi\times f}(m)m^{-s},
\end{equation}
with
\begin{equation}
  \label{eq:46}
  a_{\phi\times f}(p^\nu)\ =\ \sum_{j=1}^4 \tau_p(j)^\nu.
\end{equation}
The archimedean (gamma) factor is
\begin{eqnarray}
  \label{eq:41}
  L_\infty (s,\phi\times f) & \ := \ & \Gamma_\R(s+it_\phi+{\textstyle\frac{k-1}2})
  \Gamma_\R(s-it_\phi+{\textstyle\frac{k-1}2})\nonumber\\
& & \ \times \
  \Gamma_\R(s+it_\phi+{\textstyle\frac{k+1}2})
  \Gamma_\R(s-it_\phi+{\textstyle\frac{k+1}2}),\ \ \ \
\end{eqnarray}
and the completed $L$-function
\begin{equation}
  \label{eq:42}
  \Lambda(s,\phi\times f)\ := \ L_\infty(s,\phi\times f)L(s,\phi\times f)
\end{equation}
satisfies the functional equation
\begin{equation}
  \label{eq:43}
  \Lambda(s,\phi\times f) \ = \  \Lambda(1-s,\phi\times f).
\end{equation}
Note the functional equation is even. We define the archimedean parameters
$\mu_j$, $1\leq j\leq4$, to be the numbers
\begin{equation}
  \label{eq:44}
  \frac{k\pm1}2 \pm it_\phi.
\end{equation}

\subsection{Explicit Formula} \label{sec:explicit-formula-GL4} As
in \S\ref{sec:explicit-formula}, let $R>0$ be a parameter; later we
take $R = k^4$. Assuming GRH we may write the non-trivial zeros of
$L(s,\phi\times f)$ as $\rho_j=\frac12+i\gamma_j$, $j\in\Z-\{0\}$
(since all signs are even). Then
\begin{equation}
  \label{eq:47}
  \sum_j g\left(\frac{\gamma_j}{2\pi}\log R\right) \ = \
  \frac A{\log R} -
  2\sum_{p}\sum_{\nu=1}^\infty\hat g\left(\frac{\nu\log p}{\log R}\right)
  \frac{a_{\phi\times f}(p^\nu)\log p}{p^{\nu/2}\log R},
\end{equation}
where
\begin{equation}
  \label{eq:48}
  A \ = \  \int_{-\infty}^{\infty} \sum_{j=1}^4 \left(\frac{\Gamma'_\R}{\Gamma_\R}
  \left(\mu_j+\frac12 + \frac{2\pi ix}{\log R}\right) + \frac{\Gamma'_\R}{\Gamma_\R}
  \left(\overline{\mu_j}+\frac12 + \frac{2\pi ix}{\log R}\right)\right)g(x)\d x.
\end{equation}

An analogous calculation as in \S\ref{subsecgammasix} gives, up to
lower order terms, that the conductor is $\left(k^2/4\right)^2$. As
we only care about the logarithm of the conductor, we take $R =
k^4$. The contribution to the $1$-level density will be
$\widehat{g}(0)$ plus lower order terms.

\subsection{Relation of $a_{\phi\times f}$ to $\lambda_\phi$ and
  $\lambda_f$.}
\label{sec:relat-a-lambda} We consider the local parameters
$\alpha_p^{\pm1}$ at any prime~$p$ for $f$, as well as
$\beta_p^{\pm1}$ for $\phi$. The local parameters $\tau_p(j)$
($j=1,2,3,4$) for the automorphic representation associated to $L(s,
\phi\times f)$ are the four numbers $\alpha_p^{\pm1}\beta_p^{\pm1}$.
A calculation similar to but simpler than that in
\S\ref{sec:a_phi-expansion} gives
\begin{align}
  \label{eq:49}
  a_{\phi\times f}(p) &\ = \  \lambda_\phi(p)\lambda_f(p) \\
  \label{eq:50}
  a_{\phi\times f}(p^\nu) &\ = \
  (\lambda_\phi(p^\nu)-\lambda_\phi(p^{\nu-2})
  (\lambda_f(p^\nu)-\lambda_f(p^{\nu-2})), \ \nu \ge 2.
\end{align}
In particular
\begin{equation}
  \label{eq:51}
  a_{\phi\times f}(p^2) \ = \  (\lambda_\phi(p^2)-1)(\lambda_f(p^2)-1).
\end{equation}

\subsubsection{$\nu \geq 3$ Terms}

We show there is no contribution to the $1$-level density from terms
with $\nu \geq 3$ in \eqref{eq:47}. The Satake parameters are the
four numbers $\ga_p^{\pm 1}\gb_p^{\pm 1}$, each of which is bounded
by $p^{\delta}$ with $\delta < \frac{1}{6}$ by Kim-Sarnak \cite{K}.
Thus \be \frac{a_{\phi,f}(p^\nu)}{p^{\frac{\nu}{2}}} \ \ll \
p^{\nu(\delta - \foh)}, \ee and as $\delta < \frac{1}{6}$, summing
over $p$ and $\nu \geq 3$ is $O(1)$. Dividing by $\log R = \log
k^4$, we see there is no contribution.

\subsubsection{$\nu = 1$ Terms}

As $a_{\phi \times f} = \gl_\phi(p)\gl_f(p)$ and $\log R = \log
k^4$, we have \bea  -2\sum_{p}\hat g\left(\frac{\log p}{\log
R}\right) \frac{\gl_\phi(p)\gl_f(p)\log p}{\sqrt{p}\log R}. \eea
As $1 = \lambda_{f}(1)$, summing over $f\in \hkn$ yields no
contribution for small support, as $\glf(p)\glf(1)$ is a
non-diagonal term.

\subsubsection{$\nu = 2$} Terms

As $a_{\phi\times f}(p^2) = (\lambda_\phi(p^2)-1)
(\lambda_f(p^2)-1)$, we have \be -2\sum_{p}\hat g\left(\frac{2\log
p}{\log R}\right) \frac{(\lambda_\phi(p^2)-1)
(\lambda_f(p^2)-1)\log p}{p\log R}. \ee There are four types of
terms: \be \glp(p^2)\cdot\glf(p^2)\glf(1),\ \ -\glp(p^2), \ \
-\glf(p^2)\glf(1), \ \ (-1)(-1). \ee The first and third are
non-diagonal, and by the Petersson formula will not contribute for
small support. The second is diagonal; however, by GRH for
$L(s,\sym^2 \phi)$ (see \eqref{eqglsixsymphi}), this term is
$O(\frac{1}{\log R})$. We are left with the fourth piece, \be
-(-1)^2 \cdot 2\sum_{p}\hat g\left(\frac{2\log p}{\log R}\right)
\frac{\log p}{p\log R}. \ee By Lemma \ref{thmprimesums}, the
$p$-sum is $\frac{g(0)}{4}$, thus, the $\nu=2$ terms contribute,
for small support, $-\foh g(0)$.

\subsubsection{Summary}

The previous subsections proved Theorem \ref{thm:sympfam}, that as
$k\to\infty$ the $1$-level density of $\mathcal{F}_{\phi\times
H_k}$ agrees only with symplectic. We took two orthogonal families
(when $k \equiv 2 \bmod 4$ then $\hkn$ has $\soo$ symmetry, and
when $k \equiv 0 \bmod 4$ then $\hkn$ has $\soe$ symmetry), and
showed that their twists by a fixed even, full level Maass form
give a symplectic family. This should be compared to Theorem
\ref{thm:soefam}, where we twisted a symplectic family and
obtained an $\soe$ family.

\begin{rek}
The reason for the symmetry flipping can be found in
\eqref{eq:51}. For the $\nu = 2$ terms, the Maass form introduces
an extra factor of $-1$ in the diagonal contribution. This changes
the sign of the contribution from the $\nu = 2$ terms, and
switches us from symplectic to orthogonal symmetries (if we have
orthogonal symmetries, we need to evaluate the 2-level density to
determine which one as our supports are too small to distinguish
$\soe$, $\so$ and $\soo$).
\end{rek}


\section{Conclusion}
\label{sec:summary-all}

We investigated the distribution of low lying zeros for two
families. The first is a $\glsix$ family, $\{\phi \times \sym^2 f: f
\in H_k\}$; here $\phi$ is a fixed Hecke-Maass cusp form and $k \to
\infty$. Though this is a $\glsix$ family, only four of the six
Gamma factors depend on $k$, and the analytic conductor is $k^4$.
Since all elements of this family have even functional equation and
there is no natural complementary family with odd sign, a folklore
conjecture predicted that the underlying group symmetry should be
symplectic. However, the symmetry type is $\soe$, proving that low
lying zeros is more than just a theory of signs of functional
equations.

We calculated the $1$-level density for test functions $g$ such that
${\rm supp}(\widehat{g})$ is small (${\rm supp}(\widehat{g}) \subset
(-\frac5{24},\frac5{24}$)). For such small support, only the
diagonal terms in the Petersson formulas contribute. Thus we can
eliminate two of the five classical compact groups, namely
symplectic and unitary. Unfortunately, since the support is
significantly less than $(-1,1)$, all three orthogonal candidates
are still possible; however, as all members of the family have even
functional equation, we do not expect the underlying symmetry to be
either $\so$ or $\soo$. Observe that the reason for the flipping of
symmetry from symplectic to (some flavor of) orthogonal is that the
contribution from the squares of primes (which is what is
responsible in any case for the term $\pm\foh g(0)$) changes sign.

More precisely, in \cite{ILS} for the family $\{\sym^2 f\}$, there
is a contribution of $-\foh \hg{0}$, arising from the diagonal
term $+1$ in $\lambda_f(p^4)-\lambda_f(p^2)+1$. In our family,
this term is multiplied by the factor $\lambda_\phi(p^2)-1$ (see
\eqref{eq:25}), and the $-1$ results in a diagonal contribution of
opposite sign, hence the symmetry flipping.

To discard the $\so$ and $\soo$ symmetries, we calculated the
$2$-level density.  It is shown in \cite{Mil1,Mil2} that for
arbitrarily small support the $2$-level densities of the three
orthogonal groups are distinguishable. We see that our answer
agrees only with $\soe$, further supporting the claim that the
symmetry group of this family is $\soe$.

Our second example is a $\glfour$ family, $\{\phi \times f: f \in
\hkn\}$. Here, for the same reason as before, twisting flips the
symmetry (this time from orthogonal to symplectic). In a subsequent
paper, we will describe the interplay between twisting by a fixed
$\gln$ form (or family) and the symmetry type (in certain cases).
The arguments of this paper can be generalized to families
satisfying certain natural technical conditions.  It can be shown
that a natural ``family constant'' $c$ can be attached to a family
in such a manner that $c_{\mathcal{F}\times\mathcal{G}} =
c_{\mathcal{F}} \cdot c_{\mathcal{G}}$, where $c = 0$ $(1, -1)$ for
unitary (symplectic, orthogonal) symmetry.  Here
$\mathcal{F}\times\mathcal{G}$ is the family obtained by
Rankin-Selberg convolution of the $L$-functions in the families
$\mathcal{F}$ and $\mathcal{G}$. These results are similar in spirit
to the universality found by Rudnick and Sarnak \cite{RS} in the
$n$-level correlations of high zeros, and will be described in
further detail in \cite{DM}. Specifically, it again seems that the
second moment of the Satake parameters determines the answer. For
the $\glsix$ family, the main term from averaging the second moment
of the Satake parameters (see \eqref{eq:16}) is $-1$ and leads to
$\soe$ symmetry, while in the $\glfour$ family the main term (see
\eqref{eq:51}) is $+1$ and leads to symplectic symmetry.


\appendix

\section{Gamma factors and signs of functional equations}
\label{sec:appendix}

In this appendix we derive the precise forms (equations \eqref{eq:6}
and~\eqref{eq:41}) of the gamma factors for the completed
$L$-functions $L(s,\phi\times\sym^2f)$ and $L(s,\phi\times f)$, as
well as their functional equations (equations \eqref{eq:10} and
\eqref{eq:43}). In particular, we show that both functional equations
are even.

Being Hecke eigenforms of level~1, $f$ and $\phi$ can be identified
with automorphic cuspidal representations $F$ and $\Phi$ of
$\GL_2(\bfA_\bfQ)$ with trivial central character~\cite{Gel}.  The
latter are isomorphic to restricted tensor products
\begin{equation}
  \label{eq:20}
  F\ \cong\ \sideset{}{'}\prod_v F_v \qquad \Phi\ \cong\ \sideset{}{'}\prod_v \Phi_v
\end{equation}
of representations of $GL_2(\bfQ_v)$ for each place $v$ of~$\bfQ$
(here $v=p$ for $p$ prime, or $v=\infty$, in which case
$\bfQ_{\infty}=\bfR$).  In the case at hand, as every (finite) prime
$p$ is unramified, the corresponding principal series
representations of $\GL_2(\bfQ_p)$ have associated
$\SL_2(\bfC)$-conjugacy classes%
\footnote{$A^\natural$ denotes the conjugacy class of~$A$.}
\begin{equation}
  \label{eq:54}
  F_p\ \leftrightarrow\
  \begin{pmatrix}
    \alpha_p&0\\0&\alpha_p^{-1}
  \end{pmatrix}^{\natural},\qquad
  \Phi_p\ \leftrightarrow\
  \begin{pmatrix}
    \beta_p&0\\0&\beta_p^{-1}
  \end{pmatrix}^{\natural}.
\end{equation}
Denote by $M(\alpha_p),M(\beta_p)$ the matrices in~(\ref{eq:54}).

Taking the symmetric square of the standard representation of
$\GL_2(\bfR)$, one obtains the conjugacy class
$\sym^2M(\alpha_p)^\natural=\diag(\alpha_p^2,1,\alpha_p^{-2})^\natural$,
whence the Satake parameters of $\sym^2f$.  Similarly, the conjugacy
classes of $M(\alpha_p)\otimes M(\beta_p)$ in $\GL_4(\bfR)$ and
$\sym^2M(\alpha_p)\otimes M(\beta_p)$ in $\GL_6(\bfR)$ define the
Satake parameters of $\phi\times f$ and~$\phi\times\sym^2f$.  The
local $L$-factors are defined in terms of these Satake parameters in
the usual manner, and the product of all local factors defines the
(incomplete) $L$-functions~\eqref{eq:4} and~(\ref{eq:39}).


The representations $F_\infty$ and $\Phi_\infty$ are the discrete series
representation of weight $k$, and the representation
$I(|\cdot|^{it},|\cdot|^{-it})$ of $\GL_2(\bfR)$, respectively (recall that
$\frac14+t^2$ is the Laplacian eigenvalue of $\phi$.)%
\footnote{$I(|\cdot|^{it},|\cdot|^{-it})$ is the unitary induction from the
  group $Q$ of upper-triangular matrices to $\GL_2(\bfR)$ of the
  representation $ \begin{pmatrix} a&*\\&d\end{pmatrix}
  \mapsto|a|^{it}|d|^{-it}$.}  Selberg proved that $t$ is real for Maass forms
of level~$1$ (and conjectured that $t$ is still real for forms of
any weight). The first published proof is due to Roelcke \cite{Roe};
see also \cite{Iw1}.

Proofs of equations~(\ref{eq:6}) and~(\ref{eq:41}) involve
parametrizing the representations $F_\infty$ and~$\Phi_\infty$,
through the Langlands correspondence, by semisimple representations
of the Weil group $W_\bfR$ (see \cite{Kn} and \cite{CM}). We number
the discrete series as in~\cite{Kn}, (see the note at the top of
page~1588 of~\cite{CM}), so replacing $k$ by $\ell+1$ in what
follows would make our notation agree with that of~\cite{CM}.  Then
(cf., equations (3.2) and (3.3) of~\cite{Kn})
\begin{align}
  \label{eq:56}
  F_\infty&\ \leftrightarrow\ \rho_{(k-1,0)},\\
  \Phi_\infty&\ \leftrightarrow\ \rho_{(+,it)}\oplus\rho_{(+,-it)}.
\end{align}
We have denoted by $\rho_{(a,b)}$ the semisimple representation of
$W_\bfR$ with parameters $(a,b)$.  The known cases of
functoriality~\cite{GeJa,Ra,KiSh} imply
the existence of automorphic (cuspidal) representations $\sym^2F$,
$\Phi\times\sym^2F$, and $\Phi\times F$ such that, by the
archimedean Langlands correspondence,
\begin{align}
\label{eq:60}
  (\sym^2F)_\infty &\ \leftrightarrow\ \sym^2\rho_{(k-1,0)}\\
  (\Phi\times\sym^2F)_\infty
  &\ \leftrightarrow\
  \big(\rho_{(+,it)}\oplus\rho_{(+,-it)}\big)\otimes\sym^2\rho_{(k-1,0)}\\
  (\Phi\times F)_\infty &\ \leftrightarrow\
  \big(\rho_{(+,it)}\oplus\rho_{(+,-it)}\big)\otimes\rho_{(k-1,0)}.
\end{align}
By Proposition~3.1 of~\cite{CM},%
\footnote{Note that the weight $k$ of our $f$ is always even.  Thus, the
  appearance of $(-,0)$ in~(\ref{eq:57}) is due to the fact that
  $(-1)^{k-1}=-1$.}
\begin{equation}
  \label{eq:57}
  \sym^2\rho_{(k-1,0)}\ \cong\ \rho_{(-,0)}\oplus\rho_{(2k,0)}.
\end{equation}
Moreover, it is easily checked that
\begin{align}
  \label{eq:58}
  \rho_{(+,\pm it)}&\otimes\rho_{(-,0)}\ \cong\ \rho_{(-,\pm it)} \\
  \rho_{(+,\pm it)}&\otimes\rho_{(\ell,0)}\ \ \cong\ \rho_{(\ell,\pm it)}.
\end{align}
Hence,
\begin{align}
  \label{eq:59}
  (\Phi\times\sym^2F)_\infty &\ \leftrightarrow\
  \rho_{(-,it)}\oplus\rho_{(-,-it)}
  \oplus\rho_{(2k-2,it)}\oplus\rho_{(2k-2,-it)}\\
  \label{eq:62}
  (\Phi\times F)_\infty &\ \leftrightarrow\
  \rho_{(k-1,it)}\oplus\rho_{(k-1,-it)}.
\end{align}
The archimedean (gamma) factors can be found using these
decompositions and the local Langlands correspondence. In terms of
irreducible semisimple representations of $W_\bfR$, we have
\begin{equation}
\label{eq:63}
  L(s,\rho)\ =\
  \begin{cases}
    \Gamma_\bfR(s\pm it) & \rho\ = \ \rho_{(+,\pm it)}, \\
    \Gamma_\bfR(s\pm it+1) & \rho \ = \ \rho_{(-,\pm it)}, \\
    \Gamma_\bfR(s\pm it+\frac\ell2)\Gamma_\bfR(s\pm it+\frac\ell2+1)
    & \rho\ = \ \rho_{(\ell,\pm it)}.
  \end{cases}
\end{equation}
These local factors are multiplicative under direct sums of
representations of $W_\bfR$, so definitions (\ref{eq:6})
and~(\ref{eq:41}) are consistent with \eqref{eq:59} and~\eqref{eq:62}
via~(\ref{eq:63}).

Since all automorphic representations under discussion are
self-con\-tra\-gred\-ient, the functional equations relate each
$L$-function to itself as $s\mapsto1-s$.  In general, the root number
$\varepsilon(s,\pi)$ associated to an automorphic cuspidal
representation $\pi$ is a product
\begin{equation*}
  \varepsilon(s,\pi) = \prod_v \varepsilon(s,\pi_v)
\end{equation*}
of local root numbers.%
\footnote{We have omitted the dependence of the local root numbers on
  the choice of additive character $\psi$ of $\bfA_\bfQ$.}  For
self-contragredient representations $\pi=\tilde\pi$ the
$\varepsilon$-factor agrees with the sign of the functional equation
(up to a factor $Q^s$ which is not present in the level-$1$ case that
concerns us).  Moreover, local root numbers can be computed via the
local Langlands correspondence as root numbers associated to Weil
group representations.  Additionally, $\varepsilon(s,\pi_p)=1$ at any
prime $p$ such that $\pi_p$ is unramified.\footnote{ We assume
  $\psi_p$ is unramified at all primes $p$.}  The local root number of
an irreducible representation $\rho$ of $W_\bfR$ is
\begin{align}
  \label{eq:61}
  \varepsilon(s,\rho) =
  \begin{cases}
    1 & \rho=\rho_{(+,\pm it)}, \\
    i & \rho=\rho_{(-,\pm it)}, \\
    i^{\ell+1} & \rho=\rho_{(\ell,\pm it)}.
  \end{cases}
\end{align}
From~(\ref{eq:61}) and~\eqref{eq:59}, (\ref{eq:62}) we obtain
\begin{align}
  \varepsilon((\Phi\times\sym^2F)_\infty)
  &\ =\ i\cdot i\cdot i^{2k-1}\cdot i^{2k-1}\ =\ +1,
\end{align} and, since $k$ is even,  \begin{align} \varepsilon((\Phi\times F)_\infty)
  &\ =\ i^k\cdot i^k\ =\ +1. \end{align}

Since all (finite) primes $p$ are unramified for $\Phi\times F$ and
$\Phi\times\sym^2F$, we conclude
\begin{equation}
  \label{eq:64}
  \varepsilon(s,\Phi\times\sym^2F)\ =\ \varepsilon(s,\Phi\times F)\ =\ +1,
\end{equation}
so the global functional equations have even sign, proving~(\ref{eq:10})
and~(\ref{eq:43}).


\bigskip


\begin{thebibliography}{CFKRS} 


\bibitem[Bu1]{Bu1}
\newblock D.~Bump, \emph{Automorphic forms on $\text{GL}(3,\R)$},
Lecture Notes in Mathematics, 1083, Springer-Verlag, Berlin, 1984.

\bibitem[Bu2]{Bu2}
\newblock D.~Bump, \emph{The Rankin-Selberg method: a survey},
Number theory, trace formulas and discrete groups (Oslo, 1987),
Academic Press, Boston, MA, 1989, 49--109.

\bibitem[BDHI]{BDHI}
\newblock D.~Bump, W.~Duke, J.~Hoffstein, and H.~Iwaniec,
\emph{An estimate for the Hecke eigenvalues of Maass forms},
Int.~Math.~Res.~Not.~1992, no.~4, 75--81.


\bibitem[CM]{CM}
J. Cogdell and P. Michel,
\emph{On the complex moments of symmetric power $L$-functions at $s=1$},
Int.~Math.~Res.~Not.~2004, no.~31, 1561--1617.


\bibitem[CFKRS]{CFKRS} \newblock B.~Conrey, D.~Farmer, J.~P.~Keating,
  M.~Rubinstein and N.~Snaith, \emph{Integral Moments of
    $L$-Functions}, Proc. London Math. Soc. (3) \textbf{91} (2005),
  no. 1, 33--104.

\bibitem[Da]{Da}
\newblock H.~Davenport, \emph{Multiplicative Number Theory, $2$nd edition},
 Graduate Texts in Mathematics \textbf{74}, Springer-Verlag, New York,
 $1980$, revised by H. Montgomery.

\bibitem[DM]{DM} \newblock E.~Due\~nez and S.~J.~Miller, \emph{The
    effect of twisting families of $L$-functions on the underlying
    group symmetries}, preprint.

\bibitem[FI]{FI}
\newblock E. Fouvry and H. Iwaniec, \emph{Low-lying zeros of dihedral
$L$-functions}, Duke Math. J.  \textbf{116} (2003),  no.~2,
189--217.

\bibitem[Gel]{Gel} \newblock S.~Gelbart, \emph{Automorphic forms on ad{\`e}le
    groups}, Annals of Mathematics Studies, No.~83.  Princeton University
  Press, Princeton, N.J., 1975.

\bibitem[GeJa]{GeJa} \newblock S.~Gelbart and H.~Jacquet, \emph{A relation
    between automorphic representations of ${\rm GL}(2)$ and ${\rm GL}(3)$},
  Ann.~Sci.~École Norm.~Sup.~\textbf{4} 11 (1978), no.~4, 471--542.

\bibitem[GR]{GR}
\newblock I.~Gradshteyn and I.~Ryzhik, \emph{Tables of Integrals,
Series, and Products}, New York, Academic Press, 1965.

\bibitem[G\"u]{Gu}
\newblock A. G\"ulo\u{g}lu, \emph{Low Lying Zeros of Symmetric
Power $L$-Functions}, IMRN, Vol. 2005, Issue 9, 517--550.

\bibitem[HW]{HW}
\newblock G.~Hardy and E.~Wright, \emph{An Introduction to the
Theory of Numbers}, fifth edition, Oxford Science Publications,
Clarendon Press, Oxford, 1995.

\bibitem[Hej]{Hej}
\newblock D.~Hejhal, \emph{On the triple correlation of zeros of
the zeta function}, Internat. Math. Res. Notices $1994$, no. $7$,
294--302.

\bibitem[HR1]{HR1} \newblock C. Hughes and Z. Rudnick, \emph{Mock
    Gaussian behaviour for linear statistics of classical compact
    groups}, J. Phys. A \textbf{36} (2003), 2919--2932.

\bibitem[HR2]{HR2}
\newblock C. Hughes and Z. Rudnick, \emph{Linear Statistics of
Low-Lying Zeros of $L$-functions}, Quart. J. Math. Oxford
\textbf{54} (2003), 309--333.

\bibitem[HM]{HM}
\newblock C. Hughes and S. J. Miller, \emph{Low-lying zeros of $L$-functions
with Orthogonal symmetry}, preprint.
(http://arxiv.org/pdf/math.NT/0507450).


\bibitem[Iw1]{Iw1} \newblock H.~Iwaniec, \emph{Introduction to the Spectral
    Theory of Automorphic Forms}, Biblioteca de la Revista Matem{\'a}tica
  Iberoamericana, Madrid, 1995.

\bibitem[Iw2]{Iw2}
H. Iwaniec, \emph{Topics in Classical Automorphic Forms}, Graduate
Studies in Mathematics, Vol. 17, AMS, Providence, RI, 1997.

\bibitem[ILS]{ILS}
\newblock H.~Iwaniec, W.~Luo and P.~Sarnak, \emph{Low lying zeros of
families of $L$-functions}, Inst. Hautes \'Etudes Sci. Publ. Math.
\textbf{91} (2000), 55--131.

\bibitem[KS1]{KS1}
\newblock N.~Katz and P.~Sarnak, \emph{Random Matrices, Frobenius
Eigenvalues and Monodromy}, AMS Colloquium Publications
\textbf{45}, AMS, Providence, $1999$.

\bibitem[KS2]{KS2}
\newblock N.~Katz and P.~Sarnak, \emph{Zeros of zeta functions and symmetries},
Bull. AMS \textbf{36} (1999), 1--26.

\bibitem[KeSn]{KeSn}
\newblock J.~P.~Keating and N.~C.~Snaith, \emph{Random
matrices and $L$-functions}, Random matrix theory, J. Phys. A
\textbf{36} (2003), no.~12, 2859--2881.

\bibitem[K]{K} \newblock H.~Kim, \emph{Functoriality for the exterior
    square of $GL_2$ and the symmetric fourth of $GL_2$}, Jour. AMS
  \textbf{16} (2003), no. 1, 139--183.

\bibitem[KiSh]{KiSh}
\newblock H.~Kim and F.~Shahidi, \emph{Functorial products for
${\rm GL}\sb 2\times{\rm GL}\sb 3$ and the symmetric cube for ${\rm
GL}\sb 2$}, Annals of Math. \textbf{155} (2002), 837--893.

\bibitem[Kn]{Kn} \newblock A.~W.~Knapp, \emph{Local Langlands Correspondence:
    The Archimedean Case,} Proc.~Symp. Pure Math.~\textbf{55} (1994), Part~2,
  393--410.

\bibitem[LS]{LS} \newblock W.~Luo and P.~Sarnak, \emph{Mass equidistribution
    for Hecke eigenforms}, Comm. Pure Appl. Math. \textbf{56} (2003), no.~7,
  874--891.

\bibitem[Mil1]{Mil1} \newblock S.~J.~Miller, \emph{$1$- and $2$-Level
    Densities for Families of Elliptic Curves: Evidence for the Underlying
    Group Symmetries}, Ph.~D.~Thesis, Princeton University, $2002$.
    \\
  (http://www.math.princeton.edu/$\sim$sjmiller/thesis/thesis.html).

\bibitem[Mil2]{Mil2} \newblock S. J. Miller, \emph{$1$- and $2$-Level
    Densities for Families of Elliptic Curves: Evidence for the Underlying
    Group Symmetries}, Compositio Mathematica \textbf{104} (2004), 952--992.


\bibitem[Mon]{Mon}
\newblock H.~Montgomery, \emph{The pair correlation of zeros of the zeta
function}, Analytic Number Theory, Proc. Sympos. Pure Math.
\textbf{24}, Amer. Math. Soc., Providence, 1973, 181--193.


\bibitem[Od]{Od} \newblock A.~Odlyzko, \emph{The $10^{22}$-nd zero of the
    Riemann zeta function}, Proc.  Conference on Dynamical, Spectral and
  Arithmetic Zeta-Functions, M. van Frankenhuysen and M. L. Lapidus, eds.,
  Amer. Math. Soc., Contemporary Math. series, $2001$.\\
  (http://www.research.att.com/$\sim$amo/doc/zeta.html).

\bibitem[Ra]{Ra} \newblock D.~Ramakrishnan, \emph{Modularity of the
  Rankin-Selberg $L$-series, and multiplicity one for ${\rm SL}(2)$,}
  Ann.~of~Math. (2) \textbf{152} (2000), no.~1, 45--111.

\bibitem[Roe]{Roe}
\newblock W. Roelcke, \emph{$\ddot{U}$ber die Wellengleichung be
Grenzkreisgruppen erster Art.}, S.-B Heidelberger Akad. Wiss. Math.
Nat. Kl., 1956, 4 Abh.

\bibitem[Ro]{Ro}
\newblock E. Royer, \emph{Petits z\'{e}ros de fonctions $L$
de formes modulaires}, Acta Arith. \textbf{99} (2001),  no. 2,
147--172.

\bibitem[Ru1]{Ru1}
\newblock M.~Rubinstein, \emph{Evidence for a spectral
interpretation of the zeros of $L$-functions}, Ph.~D. Thesis,
Princeton University, $1998$.\\
(http://www.math.uwaterloo.ca/$\sim$mrubinst/thesis/thesis.html).

\bibitem[Ru2]{Ru2}
\newblock M. Rubinstein, \emph{Low-lying zeros of $L$-functions and
random matrix theory}, Duke Math. J. \textbf{109} (2001), no.~1,
147--181.

\bibitem[RS]{RS}
\newblock Z.~Rudnick and P.~Sarnak, \emph{Zeros of principal $L$-functions
 and random matrix theory}, Duke Journal of Math. \textbf{81}
(1996), 269--322.

\bibitem[Si]{Si}
\newblock J. Silverman, \emph{The average rank of an algebraic
family of elliptic curves}, J. reine angew. Math. \textbf{504}
(1998), 227--236.


\bibitem[We]{We}
\newblock A.~Weil, \emph{Sur les `formules explicites' de la
th\'eorie des nombres premiers}, Comm. S\'em. Math. Univ. Lund (Medd.
Lunds Univ. Mat. Sem.) 1952, Tome Supplementaire, 252--265.

\bibitem[Yo]{Yo}
\newblock M. Young, \emph{Low Lying Zeros of Families of Elliptic
Curves}, J. Amer. Math. Soc. \textbf{19} (2006), no. 1, 205--250.



\end{thebibliography}
\end{document}